# Directed homotopy theory, I. The fundamental category [*]

Marco Grandis

*Dipartimento di Matematica, Università di Genova, Via Dodecaneso 35, 16146-Genova, Italy.*
*e-mail: grandis@dima.unige.it*

**Abstract**. Directed Algebraic Topology is beginning to emerge from various applications.

The basic structure we shall use for such a theory, a *d-space*, is a topological space equipped with a family of *directed paths*, closed under some operations. This allows for *directed homotopies*, generally non reversible, represented by a cylinder and cocylinder functors. The existence of 'pastings' (colimits) yields a geometric realisation of cubical sets as d-spaces, together with homotopy constructs which will be developed in a sequel. Here, the *fundamental category* of a d-space is introduced and a 'Seifert - van Kampen' theorem proved; its homotopy invariance rests on 'directed homotopy' of categories. In the process, new shapes appear, for d-spaces but also for small categories, their elementary algebraic model.

Applications of such tools are briefly considered or suggested, for objects which model a directed image, or a portion of space-time, or a concurrent process.

**Introduction**

Directed Algebraic Topology is a recent subject, for which some references are given below. Its domain should be distinguished from classical Algebraic Topology by the principle that *directed spaces have privileged directions and directed paths therein need not be reversible*. Its homotopical tools, corresponding to ordinary homotopies, fundamental group and fundamental n-groupoids, should be similarly 'non-reversible': *directed homotopies, fundamental monoids and fundamental n-categories*. Its applications will deal with domains where privileged directions appear, like concurrent processes, traffic networks, space-time models, etc. Formally, new 'shapes' arise, with an aesthetic appeal often strengthened by the logical necessity of homotopy constructs (cf. 1.2, 1.6, 4.4).

As an elementary example of the notions and applications we are going to treat, consider the following (compact) zones X, Y of the plane, equipped with the order $\leq_B$. We shall see that there are, respectively, 3 or 4 homotopy classes of 'directed paths' from a to b, in the fundamental categories $\uparrow\Pi_1(X)$, $\uparrow\Pi_1(Y)$, while there are none from b to a, and every loop is constant (the prefixes $\uparrow$, d- are used to distinguish a directed notion from the corresponding 'reversible' one)

(1)  $(x, y) \leq_B (x', y') \quad \Leftrightarrow \quad |y' - y| \leq x' - x,$

[*] Work supported by MIUR Research Projects.

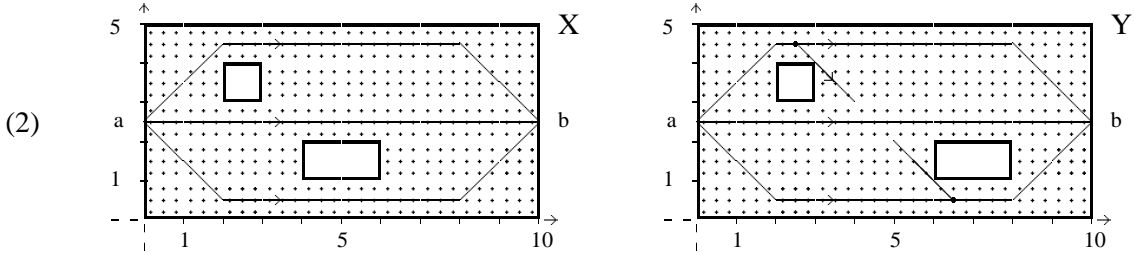

(2)

   First, these 'directed spaces' can be viewed as representing a stream with two islands; the order expresses the fact that lateral movements have an *upper bound for velocity*, equal to the speed of the stream. Secondly, one can view the abscissa as time, the ordinate as position in a 1-dimensional physical medium and the order as the possibility of going from (x, y) to (x', y') with velocity ≤ 1 (with respect to a 'rest frame', linked with the medium). The two forbidden rectangles are now linear obstacles in the medium, with a limited duration in time. Finally, our figures can be viewed as execution paths of concurrent automata, as in [FGR], fig. 14. In all these cases, the fundamental category distinguishes between obstructions (islands, temporary obstacles, conflict of resources) which intervene together (at the left) or one after the other (at the right). Note that, even if other cases can exhibit non trivial loops and *vortices* (1.6), the fundamental monoids $\uparrow\pi_1(X, x) = \uparrow\Pi_1(X)(x, x)$ often carry a very minor part of the information of $\uparrow\Pi_1(X)$.

   Now, to develop the basic theory of directed homotopy, corresponding to the ordinary theory in **Top** (topological spaces), we have to choose a precise notion of 'directed space'. We will use, in this sense, a topological space X equipped with a set dX of *directed paths* [0, 1] → X, closed under: constant paths, increasing reparametrisation and concatenation. Such objects, called *directed spaces* or *d-spaces*, form a category d**Top** which has general properties similar to **Top**: limits and colimits exist and are easily computed (1.1), and there is an *exponentiable* directed interval (1.7, 2.2).

   Direction is quite different from *orientation* (1.3a); the links with ordering are more subtle. Various d-spaces of interest derive from an ordinary space equipped with an order relation (1.4a), as in the case of the *directed interval* $\uparrow\mathbf{I} = \uparrow[0, 1]$; or, more generally, from a space equipped with a *local preorder* (1.4b), as for the directed circle $\uparrow\mathbf{S}^1$. However, the frame of locally preordered spaces would be insufficient for our purposes: they do not have general 'pastings' (colimits, cf. 4.6), and it would not be possible to form there the homotopy constructs of Part II: homotopy pushouts, mapping cones and suspension (1.4b, 1.6). Another interesting directed structure, Kelly's *bitopological spaces* [Ke], lacks path-objects (1.4c).

   The theory developed here is essentially based on the standard directed interval $\uparrow\mathbf{I}$, the (directed) *cylinder* $\uparrow I(X) = X \times \uparrow\mathbf{I}$ and its right adjoint, the *path functor* $\uparrow P(X) = X^{\uparrow\mathbf{I}}$ (2.1-2). Such functors, with a structure consisting of faces, degeneracy, connections and interchange, satisfy the axioms of an IP-*homotopical category*, as studied in [G2] for a different case of directed homotopy, cochain algebras; moreover, here, paths and homotopies can be concatenated. The theory produces two congruences on d**Top**, *d-homotopy* and *reversible d-homotopy* (2.3-4). The fundamental monoid and the functors $\uparrow\Pi_1(X)(x, x')$ are strictly invariant up to (bi)pointed d-homotopy (3.3). The invariance of the fundamental category as a whole, proved in Theorem 3.2, is more delicate, being based on parallel notions of *d-homotopy* and *reversible d-homotopy* for categories (4.1): the latter amounts to ordinary equivalence but the former is coarser (and gives a classification of categories which might be of interest in itself). See also the comments in 3.5.

Section 1 begins with basic properties and examples of d-spaces; Theorem 1.7 on exponentiable objects allows us to put a directed structure on the path-space. Section 2 deals with the cylinder and path functors, (directed) homotopies and double homotopies. Then, in Section 3, the fundamental category $\uparrow\Pi_1(X)$ of a d-space is defined; computations are essentially based on a van Kampen-type theorem (3.6). We end by treating, in Section 4, d-homotopy of categories, the geometric realisation of a cubical set as a d-space (4.5-6) and *directed metrisability* (4.7) with respect to asymmetric distances in Lawvere's sense [La].

Notions of directed Algebraic Topology, including directed paths and homotopies, have recently appeared within the analysis of concurrent processes; such notions have been developed for classical *combinatorial structures*, like simplicial and cubical sets, for *topological spaces with local orders* and for *Chu-spaces* [FGR, Pr, Ga, GG, G5]. Higher fundamental n-categories $\uparrow\Pi_n(X)$ have been developed by the author for simplicial sets [G4]. On a more formal level, it can be noted that a setting based on the (co)cylinder functor can be effectively adapted to a situation where reversion is missing, as already showed in [G2]. Kamps-Porter's text [KP] is a general reference for such settings, which go back to Kan [Ka]; while Quillen model structures [Qn] might be less suited for the directed case.

Category theory will be used at an elementary level. Some basic facts are repeatedly used: all (categorical) *limits* (extending cartesian products and projective limits) can be constructed from products and equalisers; dually, all *colimits* (extending sums and injective limits) can be constructed from sums and coequalisers; *left adjoint* functors preserve all the existing colimits, while right adjoints preserve limits (see [Ma, Bo]). $F \dashv G$ means that $F$ is left adjoint to $G$.

A *precedence* is a reflexive relation; a *preorder* is also transitive; an *order* is also anti-symmetric (and need not be total). A mapping which preserves precedence relations is said to be *increasing* (always used in the weak sense). A *map* between topological spaces is a continuous mapping. $\mathbf{I}$ is the standard euclidean interval $[0, 1]$. The letter $\alpha$ takes values $0, 1$, written $-, +$ in superscripts.

## 1. Directed topological spaces

Directed spaces are introduced, with their basic properties, their relations with other directed structures (1.4) and a first analysis of directed paths (1.5-7). Standard models are defined in 1.2.

**1.1. Basics.** A *directed topological space*, or *d-space* $X = (X, dX)$, will be a topological space equipped with a set $dX$ of (continuous) maps $a: \mathbf{I} \to X$, called *directed paths* or *d-paths*, satisfying three axioms:

(i) (*constant paths*) every constant map $\mathbf{I} \to X$ is directed,

(ii) (*reparametrisation*) $dX$ is closed under composition with increasing maps $\mathbf{I} \to \mathbf{I}$,

(iii) (*concatenation*) $dX$ is closed under path-concatenation: if the d-paths $a, b$ are consecutive in $X$ $(a(1) = b(0))$, then their ordinary concatenation $a+b$ is also a d-path

(1)   $(a+b)(t) = a(2t)$, if $0 \leq t \leq 1/2$,           $(a+b)(t) = b(2t-1)$, if $1/2 \leq t \leq 1$.

By reparametrisation, directed paths are also closed under the n-*ary concatenation* $a_1 + ... + a_n$ of consecutive paths, based on the partition $0 < 1/n < 2/n < ... < 1$.

A *directed map*, or *d-map* f: X → Y, is a continuous mapping between d-spaces which preserves the directed paths: if a ∈ dX, then fa ∈ dY. This category will be denoted as d**Top** (or ↑**Top**).

The d-structures on a space X are closed under arbitrary intersections in $\mathcal{P}(\mathbf{Top}(\mathbf{I}, X))$ and form therefore a complete lattice for the inclusion, or 'finer' relation (corresponding to the fact that idX be directed). A (directed) *subspace* X' ⊂ X has thus the restricted structure, selecting those paths in X' which are directed in X. A (directed) *quotient* X/R has the quotient structure, formed of finite concatenations of projected d-paths; in particular, for a subset A ⊂ |X|, X/A will denote the d-quotient of X which identifies all points of A. Similarly, d**Top** has all limits and colimits, constructed as in **Top** and equipped with the initial or final d-structure for the structural maps; for instance, a path $\mathbf{I} \to \Pi X_i$ is directed iff all its components $\mathbf{I} \to X_i$ are so, while a path $\mathbf{I} \to \Sigma X_i$ is directed iff it is so in some $X_i$.

The forgetful functor U: d**Top** → **Top** has adjoints $c_0 \dashv U \dashv C^0$, defined by the *d-discrete* structure of constant paths $c_0(X) = (X, |X|)$ (the finest) and, respectively, the *natural* d-structure of all paths $C^0(X) = (X, \mathbf{Top}(\mathbf{I}, X))$ (the largest). Topological spaces will generally be viewed in d**Top** via the natural embedding, which preserves products and subspaces.

Reversing d-paths, by the involution r: **I** → **I**, r(t) = 1 – t, gives the *reflected*, or *opposite*, d-space; this forms a (covariant) involutive endofunctor, called *reflection* (not to be confused with path reversion, cf. 1.5)

(2)  R: d**Top** → d**Top**,  R(X) = $X^{op}$,  (a ∈ d($X^{op}$) ⇔ $a^{op}$ = ar ∈ dX).

A d-space is *symmetric* if it is invariant under reflection. It is *reflexive*, or self-dual, if it is isomorphic to its reflection, which is more general (1.2).

**1.2. Standard models.** The euclidean spaces $\mathbf{R}^n$, $\mathbf{I}^n$, $\mathbf{S}^n$ will have their *natural* d-structure, admitting all (continuous) paths. **I** will be called the *natural* interval.

The *directed real line*, or *d-line* ↑**R**, will be the euclidean line with directed paths given by the increasing maps **I** → **R** (with respect to natural orders). Its cartesian power in d**Top**, the n-*dimensional real d-space* ↑$\mathbf{R}^n$ is similarly described (with respect to the product order, x ≤ x' iff $x_i$ ≤ $x_i'$ for all i). The *standard d-interval* ↑**I** = ↑[0, 1] has the subspace structure of the d-line; the *standard d-cube* ↑$\mathbf{I}^n$ is its n-th power, and a subspace of ↑$\mathbf{R}^n$. These d-spaces are not symmetric (for n > 0), yet reflexive; in particular, the canonical reflecting isomorphism

(1)  r: ↑**I** → R↑**I**,  t ↦ 1–t,

will play a role, in *reflecting* (not reversing!) paths and homotopies.

The *standard directed circle* ↑$\mathbf{S}^1$ will be the standard circle with the *anticlockwise structure*, where the directed paths a: **I** → $\mathbf{S}^1$ move this way, in the plane: a(t) = [1, ϑ(t)], with an increasing argument ϑ. It can also be obtained as the coequaliser in d**Top** of the following two pairs of maps

(2)  $\partial^-, \partial^+$: {∗} ⇒ ↑**I**,  $\partial^-$(∗) = 0,  $\partial^+$(∗) = 1,

(3)  id, f: ↑**R** ⇒ ↑**R**,  f(x) = x + 1.

The 'standard realisation' of the first coequaliser is the quotient ↑**I**/∂**I**, which identifies the endpoints (note that the d-quotient has the desired structure precisely because of the axiom on concatenation of d-paths). More generally, the *directed n-dimensional sphere* will be defined as the quotient of the directed cube $\uparrow\mathbf{I}^n$ modulo its (ordinary) boundary $\partial\mathbf{I}^n$, while $\uparrow\mathbf{S}^0$ has the discrete topology and the natural d-structure (obviously discrete)

(4)    $\uparrow\mathbf{S}^n = (\uparrow\mathbf{I}^n)/(\partial\mathbf{I}^n)$    (n > 0),                    $\uparrow\mathbf{S}^0 = \mathbf{S}^0 = \{-1, 1\}$.

All directed spheres are reflexive; their d-structure, further analysed in 1.6, can be described by an asymmetric distance (4.7.5). The standard circle has another d-structure of interest, induced by $\mathbf{R}\times\uparrow\mathbf{R}$ and called the *ordered circle*

(5)    $\uparrow\mathbf{O}^1 \subset \mathbf{R}\times\uparrow\mathbf{R}$,

where d-paths have to 'move up'. It is the quotient of $\uparrow\mathbf{I} + \uparrow\mathbf{I}$ which identifies lower and upper endpoints, separately. It is thus easy to guess that the *unpointed* d-suspension of $\mathbf{S}^0$ will give $\uparrow\mathbf{O}^1$, while the *pointed* one will give $\uparrow\mathbf{S}^1$, as well as all higher $\uparrow\mathbf{S}^n$ (Part II). Various versions of the projective plane will be constructed in Part II, as directed mapping cones. For the disc, see 1.6.

**1.3. Remarks.** (a) *Direction should not be confused with orientation*: every rotation of the plane preserves orientation, but only the trivial rotation preserves the directed structure of $\uparrow\mathbf{R}^2$; moreover, non-orientable surfaces like the Klein bottle have a d-structure locally isomorphic to $\uparrow\mathbf{R}^2$.

(b) A line in $\uparrow\mathbf{R}^2$ inherits the canonical d-structure (isomorphic to $\uparrow\mathbf{R}$) *if and only if* it has a positive slope (in $[0, +\infty]$); otherwise, it acquires the *discrete* d-structure (on the euclidean topology). Similarly, the d-structure induced by $\uparrow\mathbf{R}^2$ on any circle has two d-discrete arcs, where the slope is negative; $\uparrow\mathbf{S}^1$ *cannot be embedded in the directed plane* (1.6).

(c) The join of the d-structures of $\uparrow\mathbf{R}$ and $\uparrow\mathbf{R}^{op}$ is not the natural $\mathbf{R}$, but a finer structure $\mathbf{R}^\sim$: a d-path there is a *piecewise monotone* map $[0, 1] \to \mathbf{R}$, i.e. a finite concatenation of increasing and decreasing maps. The *reversible interval* $\mathbf{I}^\sim \subset \mathbf{R}^\sim$ will be of interest, for reversible paths.

(d) For a *d-topological group* G, one should require that the structural operations be directed maps $G\times G \to G$ and $G \to G^{op}$. This is the case of $\uparrow\mathbf{R}^n$ and $\uparrow\mathbf{S}^1$.

**1.4. Preorders and bitopologies.** We discuss now three other possible notions of directed topology, in increasing order of generality and linked by forgetful functors to d-spaces

(1)    $p\mathbf{Top} \subset lp\mathbf{Top} \to b\mathbf{Top} \to d\mathbf{Top}$.

(a) Firstly, a *preordered topological space* $X = (X, \leq)$ will be here a topological space equipped with a preorder relation (reflexive and transitive), under *no* coherence assumptions. Such objects, with increasing (continuous) maps, form a category p**Top** which has all limits and colimits, calculated as in **Top**, with the adequate, obvious preorder. But one cannot realise thus the directed circle and we have to look for a more general notion, localising the transitive property of preorders.

(b) A *locally preordered topological space*, or *lp-space* $X = (X, \prec)$, will have a *precedence* relation $\prec$ (reflexive) which is *locally transitive*, i.e. transitive on a suitable neighbourhood of each point; (a similar stronger notion is used in [FGR, GG] and called a *local order*: the space is equipped with an open covering and a coherent system of closed orders on such open subsets). A *map* f: $X \to Y$ is required to be *locally increasing*, i.e. to preserve $\prec$ on some neighbourhood of every $x \in X$. (Note that on a given space, infinitely many local preorders may give *equivalent* lp-structures, isomorphic via the identity. This is a minor problem: it can be mended, replacing the local preorder by its *germ*, the equivalence class for the previous relation, in the same way as a manifold structure is often defined as an equivalence class of atlases; or it can be ignored, since our mending would just replace lp**Top** with an equivalent category.)

This category lp**Top** has obvious limits and sums, but not all colimits (4.6). The forgetful functor lp**Top** → d**Top** is defined by all locally increasing paths a: $[0, 1]$ → X on the ordered interval. By compactness of $[0, 1]$ and local transitivity of $\prec$, this amounts to a continuous mapping, preserving precedence on each subinterval $[t_{i-1}, t_i]$ of a suitable decomposition $0 = t_0 < t_1 < ... < t_n = 1$. The reflection R of lp**Top** reverses the precedence relation.

A d-space will be said to be *of (pre)order type* or *of local (pre)order type* if it can be obtained, as above, from a topological space with such a structure. Thus, ↑$\mathbf{R}^n$, ↑$\mathbf{I}^n$ and the ordered circle ↑$\mathbf{O}^1$ are of order type; $\mathbf{R}^n$, $\mathbf{I}^n$ and $\mathbf{S}^n$ are of *chaotic preorder* type; the product $\mathbf{R}\times\uparrow\mathbf{R}$ is of preorder type. The d-space ↑$\mathbf{S}^1$ is of local order type, deriving from *an* anticlockwise precedence relation x $\prec_\varepsilon$ x' given by $\delta(x, x') \leq \varepsilon$, with $0 < \varepsilon < 2\pi$ (the quasi-pseudo-metric $\delta(x, x')$ measures the length of the anticlockwise arc from x to x'; see 4.7).

The higher directed spheres ↑$\mathbf{S}^n$ are *not* of local preorder type, for $n \geq 2$, essentially because lp-spaces can*not* have pointlike vortices (cf. 1.6). This geometric fact can be an advantage in other contexts (for instance, all *pre*-cubical sets – with faces but *no* degeneracies – can be realised as locally ordered spaces [FGR]). But it is at the origin of the defect of colimits recalled in the Introduction and proved below (4.6), which makes them insufficient for our purposes.

(c) Finally, a *bitopological space* (a notion introduced by J.C. Kelly [Ke] and still studied) is a set equipped with a pair of topologies $X = (X, \tau^-, \tau^+)$. Their category b**Top**, with the obvious maps – continuous with respect to *past* ($\tau^-$) and *future* topologies ($\tau^+$), separately – has all limits and colimits, calculated separately on both sorts. Reflection exchanges past and future.

The forgetful functors lp**Top** → b**Top** → d**Top** are easily defined. Given an lp-space X, a fundamental system of past or future neighbourhoods at $x_0$ derives from any fundamental system $\mathcal{V}$ of the original topology, setting

(2)   $V^- = \{x \in V \mid x \prec x_0\}$,                $V^+ = \{x \in V \mid x_0 \prec x\}$               ($V \in \mathcal{V}$).

↑**I** inherits thus the left- and right-euclidean topologies. Then, a bitopological space has a canonical d-structure, with $dX = b\mathbf{Top}(\uparrow\mathbf{I}, X)$.

Problems for establishing directed homotopy in b**Top** derive from pathologies of (say) left-euclidean topologies; in fact, for a fixed Hausdorff space A, the product $-\times A$ preserves quotients (if and) only if A is locally compact ([Mi], Thm. 2.1 and footnote (5)). Thus, the cylinder endofunctor $-\times\uparrow\mathbf{I}$ in b**Top** does not preserve colimits and has no right adjoint: the path-object is missing (and homotopy pullbacks as well, while homotopy pushouts have poor properties; cf. Part II).

**1.5. Directed paths.** A *path* in a d-space X will be a d-map a: ↑**I** → X defined on the standard d-interval. Plainly, this is the same as a structural map $a \in dX$, and will also be called a *directed* path when we want to stress the difference from ordinary paths in the underlying space UX. The path a has two endpoints, or *faces* $\partial^-(a) = a(0)$, $\partial^+(a) = a(1)$. Every point $x \in X$ has a *degenerate path* $0_x$, constant at x. A *loop* ($\partial^-(a) = \partial^+(a)$) amounts to a d-map ↑$\mathbf{S}^1$ → X (by 1.2.2).

By the very definition of d-structure (1.1), we already know that the concatenation a+b of two consecutive paths ($\partial^+$a = $\partial^-$b) is directed. This amounts to saying that, in d**Top** (as for spaces), the *standard concatenation pushout* – pasting two copies of the d-interval, one after the other – can be realised as ↑**I** itself

(1)
$$\begin{array}{ccc} \{*\} & \xrightarrow{\partial^+} & \uparrow\mathbf{I} \\ \partial^- \downarrow & \searrow^{k^-} & \downarrow k^- \\ \uparrow\mathbf{I} & \xrightarrow{k^+} & \uparrow\mathbf{I} \end{array}$$
$\qquad k^-(t) = t/2, \quad k^+(t) = (t+1)/2.$

The existence of a path in X from x to x' gives a *path preorder*, x $\preceq$ x' (x' is *reachable* from x). For an lp-space, the path preorder *implies* the transitive relation spanned by the precedence relation (by the characterisation of directed paths in 1.4b); it can be chaotic, as it happens in the directed circle. For a space of preorder type, the path preorder implies the given preorder and can fairly replace it (giving the same d-paths): for instance, in ↑**O**$^1$, the path order is strictly finer than the preorder induced by **R**×↑**R** and plainly more relevant.

The equivalence relation $\simeq$ spanned by $\preceq$ gives the partition of a d-space in its *path components*, and yields a functor

(2) $\quad \uparrow\Pi_0$: d**Top** → **Set**, $\qquad\qquad \uparrow\Pi_0(X) = |X|/\simeq .$

A non-empty d-space X is *path connected* if $\uparrow\Pi_0(X)$ is a point. Then, also the underlying space UX is path connected, while the converse is obviously false (cf. 1.3b). The directed spaces ↑**R**$^n$, ↑**I**$^n$, ↑**S**$^n$ are path connected (n > 0); but ↑**S**$^n$ is more strongly so, because already $\preceq$ is chaotic (every point can reach each other).

The path a will be said to be *reversible* if also $a^{op}(t) = a(1-t)$ is a directed path in X, or equivalently if a: **I**$^\sim$ → X is a d-map on the reversible interval (1.3c); plainly, such paths are closed under concatenation. (Requiring that a be directed on the natural interval **I** is a stronger condition, not closed under concatenation: the pasting, on a point, of two copies of **I** in d**Top** is not isomorphic to **I**; however, if X is of local preorder type, the two facts are equivalent, by the characterisation of d-paths in 1.4b.)

**1.6. Vortices, discs and cones.** Loosely speaking, a non-reversible path with equal (resp. different) endpoints can be viewed as *revealing a vortex* (resp. a *stream*). If X is of preorder type, any loop in X lives in a zone where the preorder is chaotic, and is reversible; if X is ordered, any loop is constant. An lp-space (not of preorder type) can have non-reversible loops, like ↑**S**$^1$.

We shall say that the d-space X has a *pointlike vortex* at x if every neighbourhood of x in X contains some non-reversible loop. It is easy to realise a directed disc having a pointlike vortex (see below), while ↑**S**$^1$ has none. In fact, *lp-spaces cannot have pointlike vortices*. (If X = (X, $\prec$) has a pointlike vortex at x, choose a neighbourhood V of x on which $\prec$ is transitive; then, any loop a: ↑**I** → V lives in a preordered space and is reversible.)

All higher directed spheres ↑**S**$^n$ = (↑**I**$^n$)/($\partial$**I**$^n$), for n ≥ 2, have a pointlike vortex at the class [0] (of the boundary points), as showed by the following non-reversible loops, 'arbitrarily small', in ↑**S**$^2$

(1) 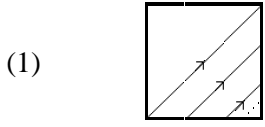

Therefore, such objects are not of local preorder type. As each other point $x \neq [0]$ has a neighbourhood isomorphic to $\uparrow\mathbf{R}^n$, this also shows that our higher d-spheres are not locally isomorphic to any fixed 'model'. (This cannot be avoided, since they are determined as pointed suspensions of $\mathbf{S}^0$.)

There are various d-structures of interest on the disc $\mathbf{B}^2$ (or compact cone). In (2), we have four cases which induce the natural structure $\mathbf{S}^1$ on the boundary: a directed path there is represented in polar coordinates as a map $[\rho(t), \vartheta(t)]$, where $\rho$ is, respectively: decreasing, increasing, constant, arbitrary. They are all of preorder type ($\rho \geq \rho'$; $\rho \leq \rho'$; $\rho = \rho'$; chaotic). One can view the first (resp. second) as a conical peak (resp. sink), and its d-structure as a decision of never going down

(2) 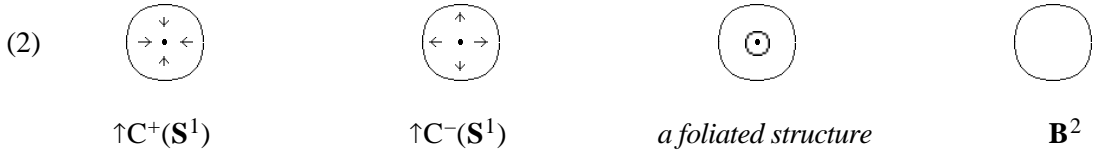

$\uparrow C^+(\mathbf{S}^1)$      $\uparrow C^-(\mathbf{S}^1)$      *a foliated structure*      $\mathbf{B}^2$

Then we have four structures which induce $\uparrow\mathbf{S}^1$ on the boundary ($\rho$ as above, $\vartheta$ increasing); all of them have a pointlike vortex at the origin (four other structures can be derived from $\uparrow\mathbf{O}^1$)

(3) 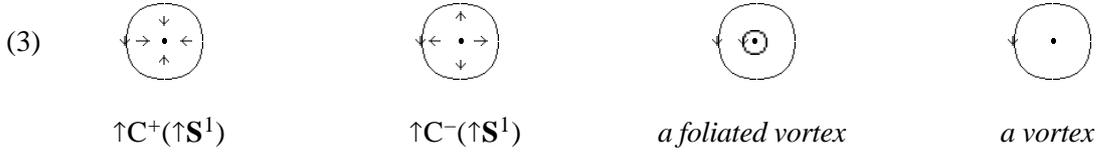

$\uparrow C^+(\uparrow\mathbf{S}^1)$      $\uparrow C^-(\uparrow\mathbf{S}^1)$      *a foliated vortex*      *a vortex*

In **Top**, the disc is the cone $C\mathbf{S}^1$, i.e. the mapping cone of $id\mathbf{S}^1$. Here, the first two cases of each row will be obtained as upper or lower directed cones (Part II), and are named accordingly.

**1.7. Theorem** (Exponentiable d-spaces). Let $\uparrow A$ be any d-structure on a locally compact Hausdorff space A. Then $\uparrow A$ is *exponentiable* in d**Top**: for every d-space Y

(1)    $Y^{\uparrow A} = d\mathbf{Top}(\uparrow A, Y) \subset \mathbf{Top}(A, UY)$,

is the set of directed maps, with the compact-open topology restricted from $(UY)^A$ and the d-structure where a path $c: \mathbf{I} \to U(Y^{\uparrow A}) \subset (UY)^A$ is directed if and only if the corresponding map $\check{c}: \mathbf{I} \times A \to UY$ is a d-map $\uparrow\mathbf{I} \times \uparrow A \to Y$.

**Proof.** It is well known that a locally compact Hausdorff space A is *exponentiable* in **Top**: the functor $-\times A: \mathbf{Top} \to \mathbf{Top}$ has a right adjoint $(-)^A: \mathbf{Top} \to \mathbf{Top}$. The space $Y^A$ is the set of maps **Top**(A, Y) with the compact-open topology; the adjunction consists of the natural bijection

(2)    $\mathbf{Top}(X, Y^A) \to \mathbf{Top}(X \times A, Y),$      $f \mapsto \check{f},$      $\check{f}(x, a) = f(x)(a).$

Now, the structure of $Y^{\uparrow A}$ defined above is well formed, as required in 1.1.i-iii.

(i) *Constant paths*. If $c: \mathbf{I} \to Y^{\uparrow A}$ is constant at the d-map $g: \uparrow A \to Y$, then $\check{c}$ can be factored as $\uparrow\mathbf{I} \times \uparrow A \to \uparrow A \to Y,$ and is directed as well.

(ii) *Reparametrisation*. For any $h: \uparrow\mathbf{I} \to \uparrow\mathbf{I}$, the map $(ch)\check{} = \check{c}.(h\times\uparrow A)$ is directed.

(iii) *Concatenation*. Let $c = c_1 + c_2: \mathbf{I} \to U(Y^{\uparrow A})$, with $\check{c}_i: \uparrow\mathbf{I}\times\uparrow A \to Y$. By the Lemma below, the product $-\times\uparrow A$ preserves the concatenation pushout 1.5.1. Therefore $\check{c}$, as the pasting of $\check{c}_i$ on this pushout, is a directed map.

Finally, we must prove that (2) restricts to a bijection between $d\mathbf{Top}(X, Y^{\uparrow A})$ and $d\mathbf{Top}(X\times\uparrow A, Y)$. In fact, we have a chain of equivalent conditions

(3)    $f: X \to Y^{\uparrow A}$ is directed,

$\forall x \in dX, \quad fx: \uparrow\mathbf{I} \to Y^{\uparrow A}$ is a d-path,

$\forall x \in dX, \quad (fx)\check{} = \check{f}.(x\times\uparrow A): \uparrow\mathbf{I}\times\uparrow A \to Y$ is directed,

$\forall x \in dX, \forall h \in d\uparrow\mathbf{I}, \forall a \in d\uparrow A, \quad \check{f}.(xh\times a): \uparrow\mathbf{I}\times\uparrow\mathbf{I} \to Y$ is directed,

$\check{f}: X\times\uparrow A \to Y$ is directed.

**1.8. Lemma.** For every d-space $\uparrow X$, the functor $\uparrow X\times -: d\mathbf{Top} \to d\mathbf{Top}$ preserves the standard concatenation pushout (1.5.1).

**Proof.** In **Top**, this is true because $X\times[0, 1/2]$ and $X\times[1/2, 1]$ form a finite closed covering of $X\times\mathbf{I}$, so that each mapping defined on the latter and continuous on such closed parts is continuous.

Consider then a map $f: X\times\mathbf{I} \to UY$ deriving from the pasting of two maps $f_0, f_1$ on the topological pushout $X\times\mathbf{I}$

(1)    $f(x, t) = f_0(x, 2t)$, for $0 \leq t \leq 1/2$,        $f(x, t) = f_1(x, 2t-1)$, for $1/2 \leq t \leq 1$.

Let now $(a, h): \uparrow\mathbf{I} \to \uparrow X\times\uparrow\mathbf{I}$ be any directed map. If the image of $h$ is contained in one half of $\mathbf{I}$, then $f.(a, h)$ is certainly directed. Otherwise, since $h$ is increasing, there is some partition $0 < t_1 < 1$ sent by $h$ to $0 < 1/2 < 1$; and we can assume that $t_1 = 1/2$ (up to composing with an automorphism of $\uparrow\mathbf{I}$). Now, the path $f.(a, h): \mathbf{I} \to UY$ is directed in $Y$, because it is the concatenation of the following two directed paths $c_i: \uparrow\mathbf{I} \to Y$

(2)    $c_1(t) = f(a(t/2), h(t/2)) = f_0(a(t/2), 2h(t/2))$,

      $c_2(t) = f(a((t+1)/2), h((t+1)/2)) = f_1(a((t+1)/2), 2h((t+1)/2) - 1)$.

## 2. Directed homotopies

The directed interval $\uparrow\mathbf{I}$ produces two adjoint endofunctors, cylinder and cocylinder, which define homotopies in $d\mathbf{Top}$. The letter $\alpha$ denotes an element of the set $\{0, 1\}$, written $-, +$ in superscripts.

**2.1. The cylinder.** The directed interval $\uparrow\mathbf{I} = \uparrow[0, 1]$ is a lattice *in* $d\mathbf{Top}$; its structure consists of two *faces* $(\partial^-, \partial^+)$, a *degeneracy* $(e)$, two *connections* or main operations $(g^-, g^+)$ and an *interchange* $(s)$

(1)    $\{*\} \underset{e}{\overset{\partial^\alpha}{\rightrightarrows}} \uparrow\mathbf{I} \overset{g^\alpha}{\leftarrow} \uparrow\mathbf{I}^2$            $s: \uparrow\mathbf{I}^2 \to \uparrow\mathbf{I}^2$,

$$\partial^\alpha(*) = \alpha, \qquad g^-(i, j) = i \vee j, \qquad g^+(i, j) = i \wedge j, \qquad s(i, j) = (j, i).$$

As a consequence, the (directed) *cylinder* endofunctor

(2) $\quad \uparrow\mathrm{I}: \mathbf{dTop} \to \mathbf{dTop}, \qquad\qquad \uparrow\mathrm{I}(-) = -\times\uparrow\mathbf{I},$

has natural transformations, which will be denoted by the same symbols and names

(3) $\quad 1 \underset{e}{\overset{\partial^\alpha}{\rightleftarrows}} \uparrow\mathrm{I} \overset{g^\alpha}{\Longrightarrow} \uparrow\mathrm{I}^2 \qquad\qquad s: \uparrow\mathrm{I}^2 \to \uparrow\mathrm{I}^2,$

and satisfy the axioms of a *cubical monad with interchange* [G2, G3].

Consecutive homotopies will be pasted via the *concatenation pushout* of the cylinder functor

(4) $\quad \begin{array}{ccc} X & \overset{\partial^+}{\longrightarrow} & \uparrow\mathrm{I}X \\ \partial^- \downarrow & & \downarrow k^- \\ \uparrow\mathrm{I}X & \underset{k^+}{\longrightarrow} & \uparrow\mathrm{I}X \end{array}$

obtained from the standard pushout 1.5.1, by applying the cartesian product $X \times -$ (Lemma 1.8). $k^\alpha X: \uparrow\mathrm{I}X \to \uparrow\mathrm{I}X$ are now two natural transformations. (The fact that pasting two copies of the cylinder gives back the cylinder is rather peculiar of spaces; e.g., it does not hold for chain complexes.)

The directed cylinder $\uparrow\mathrm{I}$ has no reversion but a *generalised reversion*, via the *reflection* of d-spaces (as for differential graded algebras [G2]; 2.2, 4.9)

(5) $\quad rX = X \times r: \uparrow\mathrm{I}.RX \to R.\uparrow\mathrm{I}X, \qquad\qquad (x, t) \mapsto (x, 1-t),$

$\qquad RrR.r = \mathrm{id}, \qquad\qquad Re.r = eR,$

$\qquad r.\partial^-R = R\partial^+, \qquad\qquad r.g^-R = Rg^+.r\uparrow\mathrm{I}.\uparrow\mathrm{I}r.$

**2.2. The path functor.** As a consequence of Theorem 1.7, the directed interval $\uparrow\mathbf{I}$ is exponentiable: the cylinder functor $\uparrow\mathrm{I} = -\times\uparrow\mathbf{I}$ has a right adjoint, the (directed) *path functor*, or *cocylinder* $\uparrow\mathrm{P}$. Explicitly, in this functor

(1) $\quad \uparrow\mathrm{P}: \mathbf{dTop} \to \mathbf{dTop}, \qquad\qquad \uparrow\mathrm{P}(Y) = Y^{\uparrow\mathbf{I}},$

the d-space $Y^{\uparrow\mathbf{I}}$ is the set of d-paths $\mathbf{dTop}(\uparrow\mathbf{I}, Y)$ with the compact-open topology (induced by the topological path-space $P(UY) = \mathbf{Top}(\mathbf{I}, UY)$) and the d-structure where a map

(2) $\quad c: \mathbf{I} \to \mathbf{dTop}(\uparrow\mathbf{I}, Y) \subset \mathbf{Top}(\mathbf{I}, UY),$

is directed iff, for all increasing maps $h, k: \mathbf{I} \to \mathbf{I}$, the derived path $t \mapsto c(h(t))(k(t))$ is in $dY$.

The lattice structure of $\uparrow\mathbf{I}$ in $\mathbf{dTop}$ produces – contravariantly – a dual structure on $\uparrow\mathrm{P}$ (a cubical *co*monad with interchange [G2, G3]); the derived natural transformations (*faces*, etc.) will be named and written as above, but proceed in the opposite direction and satisfy dual axioms (note that $\uparrow\mathrm{P}^2(Y) = Y^{\uparrow\mathbf{I}^2}$, by adjunction)

(3) $\quad \{*\} \underset{e}{\overset{\partial^\alpha}{\rightleftarrows}} \uparrow P \overset{g^\alpha}{\rightrightarrows} \uparrow P^2 \qquad\qquad s: \uparrow P^2 \to \uparrow P^2,$

$\qquad \partial^\alpha(a) = a(\alpha), \qquad e(x)(t) = x, \qquad g^-(a)(t, t') = a(t \vee t'), \ldots$

Again, the *concatenation pullback* (the object of pairs of consecutive paths) can be realised as $\uparrow PX$

(4) 
$$\begin{array}{ccc} \uparrow PX & \overset{k^+}{\longrightarrow} & \uparrow PX \\ k^- \downarrow & & \downarrow \partial^- \\ \uparrow PX & \underset{\partial^+}{\longrightarrow} & X \end{array}$$

(5) $\quad k^\alpha: \uparrow PX \to \uparrow PX, \qquad\qquad\qquad k^-(a)(t) = a(t/2), \quad k^+(a)(t) = a((t+1)/2).$

**2.3. Homotopies.** The category of d-spaces is thus an IP-*homotopical category* [G2]: it has adjoint functors $\uparrow I \dashv \uparrow P$, with the required structure (faces, etc.); it has pushouts (preserved by the cylinder) and pullbacks (preserved by the cocylinder); it has terminal and initial object. Therefore, all results of [G2] for such a structure apply (as for cochain algebras); moreover, *here we can concatenate paths and homotopies*.

A (directed) *homotopy* $\varphi: f \to g: X \to Y$ is defined as a d-map $\varphi: \uparrow IX = X \times \uparrow I \to Y$ whose two faces, $\partial^\pm(\varphi) = \varphi.\partial^\pm: X \to Y$ are $f$ and $g$, respectively. Equivalently, it is a map $X \to \uparrow PY = Y^{\uparrow I}$, with faces as above. A path is a homotopy between two points, $a: x \to x': \{*\} \to X$.

The category d**Top** will always be equipped with such homotopies and the operations produced by the (co)cylinder functor (for $\varphi: f \to g: X \to Y$; $u: X' \to X$; $v: Y \to Y'$; $\psi: g \to h: X \to Y$):

(a) *whisker composition* of maps and homotopies

$$v \circ \varphi \circ u: vfu \to vgu \qquad\qquad (v \circ \varphi \circ u = v.\varphi.\uparrow Iu: \uparrow IX' \to Y'),$$

(b) *trivial homotopies*: $\qquad 0_f: f \to f \qquad\qquad\qquad (0_f = fe: \uparrow IX \to Y),$

(c) *concatenation of homotopies*: $\quad \varphi + \psi: f \to h.$

(The horizontal composition of homotopies produces a double homotopy, 2.6a.) Concatenation is defined in the usual way, by means of the concatenation pushout (2.1.4)

(1) $\quad (\varphi+\psi)k^- = \varphi, \qquad\qquad (\varphi+\psi)k^+ = \psi \qquad\qquad\qquad (\partial^+\varphi = \partial^-\psi),$

$\qquad (\varphi+\psi)(x, t) = \varphi(x, 2t), \qquad$ for $0 \le t \le 1/2,$

$\qquad (\varphi+\psi)(x, t) = \psi(x, 2t-1), \qquad$ for $1/2 \le t \le 1.$

Directed homotopies cannot generally be *reversed*, but just *reflected* (as paths, 1.5)

(2) $\quad \varphi^{op}: Rg \to Rf: RX \to RY, \qquad\qquad \varphi^{op} = R\varphi.rX = (\uparrow I(RX) \to R(\uparrow IX) \to RY).$

*Reversible* d-homotopies $\varphi: X \times \uparrow I^\sim \to Y$, defined on the *reversible cylinder* (1.3c), have a similar structure, plus reversion; but they are rare.

The endofunctors $\uparrow I$ and $\uparrow P$ can be extended to homotopies, *via their interchange*: for $\varphi: f \to g: X \to Y$, let

(3) $\quad \uparrow\hat{I}(\varphi) = \uparrow I(\varphi).sX: \uparrow I^2(X) \to \uparrow I(Y), \qquad \uparrow\hat{P}(\varphi) = sY.\uparrow P(\varphi): \uparrow P(X) \to \uparrow P^2(Y),$

$\quad \partial^-(\uparrow IX)(\uparrow\hat{I}(\varphi)) = (\varphi \times \uparrow I).(X \times s).(X \times \uparrow I \times \partial^-) = (\varphi \times \uparrow I).(X \times \partial^- \times \uparrow I) = f \times \uparrow I = \uparrow I(f).$

**2.4. Directed homotopy equivalence.** Directed homotopies and reversible d-homotopies produce two congruence relations on d**Top**, which will be written $f \simeq g$ and $f \simeq_r g$. The second is harmless, yet generally too fine; the first is coarser and can destroy information of interest (cf. 3.5); this fact, however, can be restrained by 'fixing base points' (3.3).

The *d-homotopy preorder* $f \preceq g$, defined by the existence of a homotopy $f \to g$, is consistent with composition ($f \preceq g$ and $f' \preceq g'$ imply $f'f \preceq g'g$) but non-symmetric ($f \preceq g$ is equivalent to $Rg \preceq Rf$). It extends the path-preorder of points, $x \preceq x'$ (1.5). We shall write $f \simeq g$ the equivalence relation generated by $\preceq$: there is a finite sequence $f \preceq f_1 \succeq f_2 \preceq f_3 ... g$ (of d-maps between the same objects); it is a congruence of categories. As we have just seen, the functors $\uparrow I$ and $\uparrow P$ are *d-homotopy invariant*: they preserve the relations $\preceq$ and $\simeq$.

A *d-homotopy equivalence* will be a d-map $f: X \to Y$ having a *d-homotopy inverse* $g: Y \to X$, in the sense that $gf \simeq idX$, $fg \simeq idY$. Then we write $X \simeq Y$, and say that they are *d-homotopy equivalent*, or have the same *d-homotopy type*. A d-subspace $u: X \subset Y$ is a (directed) *deformation retract* of Y if there is a d-map $p: Y \to X$ such that $pu = idX$, $up \simeq idY$; and a *strong* deformation retract if one can choose p and the d-homotopies $up = h_0 \to h_1 \leftarrow h_2 ... h_n = idY$ so that all the latter are trivial on X. A d-space is *d-contractible* if it is d-homotopy equivalent to a point, or equivalently if it admits a deformation retract at a point. (Classical relations between ordinary homotopy equivalence and deformation retracts can be seen in [Wh], I.5.)

Plainly, all these relations imply the usual ones, for the underlying spaces. In fact, they are strictly stronger. As a trivial example, the d-discrete structure $c_0\mathbf{R}$ on the real line (where all d-paths are constant, 1.1) is not d-contractible. Less trivially, within path connected d-spaces, it is easy to show that $\mathbf{S}^1$, $\uparrow\mathbf{S}^1$ and $\uparrow\mathbf{O}^1$ are not d-homotopy equivalent: a directed map $\mathbf{S}^1 \to \uparrow\mathbf{O}^1$ or $\uparrow\mathbf{S}^1 \to \uparrow\mathbf{O}^1$ must stay in the left or right half of $\uparrow\mathbf{O}^1$, whence its underlying map is homotopically trivial. And a d-map $\mathbf{S}^1 \to \uparrow\mathbf{S}^1$ is necessarily constant.

Reversible d-homotopies (2.3) give a stronger congruence $f \simeq_r g$, and related notions: *reversible d-homotopy equivalence* ($X \simeq_r Y$) etc. The directed interval $\uparrow\mathbf{I}$ is d-contractible, but not reversibly so. Of the four d-structures considered in 1.6.2 (or 1.6.3) for the disc, the first and second are just d-contractible, the third is not, the last is reversibly d-contractible.

Each d-space A gives a covariant 'representable' *homotopy functor* (and a contravariant one)

(1) $\quad [A, -]: d\mathbf{Top} \to \mathbf{Set}, \qquad\qquad ([-, A]: d\mathbf{Top}^{op} \to \mathbf{Set}),$

where [A, X] denotes the set of d-homotopy classes of maps $A \to X$. These functors are plainly d-homotopy invariant: if $f \simeq g$ in d**Top**, then $[A, f] = [A, g]$. In particular, $\uparrow\Pi_0(X) = [\{*\}, X]$ (1.5); also $[\uparrow\mathbf{S}^1, -]$, $[\mathbf{S}^1, -]$, $[\uparrow\mathbf{O}^1, -]$ express invariants of interest; the first gives the set of *homotopy classes of directed* free *loops* (not to be confused with the fundamental monoid, 3.3).

**2.5. Double homotopies and 2-homotopies.** Roughly speaking, double homotopies (and double paths, in particular) behave as in **Top**, as long as we work on the ordered square $\uparrow\mathbf{I}^2$ via increasing maps. The *second order cylinder* $\uparrow I^2X = X \times \uparrow\mathbf{I}^2$ has four 1-dimensional faces, written

(1) $\partial_1^\alpha = \uparrow\!\!I\partial^\alpha = X\times\partial^\alpha\times\uparrow\!\!I : \uparrow\!\!IX \to \uparrow\!\!I^2 X,$ $\qquad \partial_1^\alpha(x, t) = (x, \alpha, t),$

$\partial_2^\alpha = \partial^\alpha \uparrow\!\!I = X\times\uparrow\!\!I\times\partial^\alpha : \uparrow\!\!IX \to \uparrow\!\!I^2 X,$ $\qquad \partial_2^\alpha(x, t) = (x, t, \alpha).$

A *double homotopy* is a map $\Phi: X\times\uparrow\!\!I^2 = \uparrow\!\!I^2 X \to Y$ (or, equivalently, $X \to \uparrow\!\!P^2 Y$); it has four faces, which will be drawn as below

(2) $\partial_1^\alpha(\Phi) = \Phi.\partial_1^\alpha = \Phi.(X\times\partial^\alpha\times\uparrow\!\!I),$ $\qquad \partial_2^\alpha(\Phi) = \Phi.\partial_2^\alpha = \Phi.(X\times\uparrow\!\!I\times\partial^\alpha),$

$$\begin{array}{ccc} & \partial_2^-\Phi & \\ f & \longrightarrow & h \\ \partial_1^-\Phi \downarrow & \Phi & \downarrow \partial_1^+\Phi \\ k & \longrightarrow & g \\ & \partial_2^+\Phi & \end{array} \qquad\qquad \begin{array}{ccc} \bullet & \longrightarrow & 1 \\ \downarrow & & \\ 2 & & \end{array}$$

and four vertices, $\partial^-\partial_1^-(\Phi) = f = \partial^-\partial_2^-(\Phi)$, etc. The concatenation, or pasting, of double homotopies *in direction* 1 or 2 is defined as usual (under the obvious boundary conditions) and satisfies a strict *middle-four interchange property*

(3) $(A +_1 B) +_2 (C +_1 D) = (A +_2 C) +_1 (B +_2 D),$

$$\begin{array}{ccccc} \bullet & \to & \bullet & \to & \bullet \\ \downarrow & A & \downarrow & B & \downarrow \\ \bullet & \to & \bullet & \to & \bullet \\ \downarrow & C & \downarrow & D & \downarrow \\ \bullet & \to & \bullet & \to & \bullet \end{array} \qquad\qquad \begin{array}{ccc} \bullet & \longrightarrow & 1 \\ \downarrow & & \\ 2 & & \end{array}$$

A (directed) *2-homotopy* $\Phi: \varphi \to \psi: f \to g: X \to Y$ is a double homotopy whose faces $\partial_1^\alpha$ are degenerate, while the faces $\partial_2^\alpha$ are $\varphi, \psi$ (the other choice being equivalent, by interchange)

(4) $$\begin{array}{ccc} & \varphi & \\ f & \longrightarrow & g \\ 0_f \downarrow & \Phi & \downarrow 0_g \\ f & \longrightarrow & g \\ & \psi & \end{array} \qquad \begin{array}{ll} \partial_2^-(\Phi) = \varphi, & \partial_2^+(\Phi) = \psi, \\ \partial^-\varphi = f = \partial^-\psi, & \partial^+\varphi = g = \partial^+\psi, \\ \partial_1^-(\Phi) = 0_f, & \partial_1^+(\Phi) = 0_g. \end{array}$$

Such particular double homotopies are closed under 1- and 2-pasting (also because $0_f + 0_f = 0_f$). The preorder $\varphi \preceq_2 \psi$ (there is a 2-homotopy $\varphi \to \psi$) spans an equivalence relation $\simeq_2$.

**2.6. Constructing double homotopies.** (a) Two 'horizontally' consecutive d-homotopies

(1) $\varphi: f^- \to f^+: X \to Y,$ $\qquad\qquad \psi: g^- \to g^+: Y \to Z,$

can be composed, to form a double homotopy $\psi \circ \varphi$

(2) $$\begin{array}{ccc} & g^-\circ\varphi & \\ g^-f^- & \longrightarrow & g^-f^+ \\ \psi\circ f^- \downarrow & \psi\circ\varphi & \downarrow \psi\circ f^+ \\ g^+f^- & \longrightarrow & g^+f^+ \\ & g^+\circ\varphi & \end{array} \qquad \begin{array}{l} \psi\circ\varphi = \psi.(\varphi\times\uparrow\!\!I): X\times\uparrow\!\!I^2 \to Y\times\uparrow\!\!I \to Z, \\ \partial_1^\alpha(\Phi) = \psi.(\varphi\partial^\alpha\times\uparrow\!\!I) = \psi\circ f^\alpha, \\ \partial_2^\alpha(\Phi) = \psi.(\varphi\times\partial^\alpha) = g^\alpha\circ\varphi. \end{array}$$

(Together with the whisker composition, in 2.3, this is a particular instance of the cubical enrichment produced by the (co)cylinder functor: composing a p-uple homotopy $\Phi: \uparrow I^p X \to Y$ with a q-uple one $\Psi: \uparrow I^q Y \to Z$ gives a (p+q)-uple homotopy $\Psi \circ \Phi = \Psi . \uparrow I^q \Phi = \uparrow I^p \Psi . \Phi$.)

(b) *Acceleration*. For every homotopy $\varphi: f \to g$, there are *acceleration* 2-homotopies

(3)   $\Theta': 0_f + \varphi \to \varphi,$  $\qquad\qquad\qquad\qquad \Theta'': \varphi \to \varphi + 0_g,$

(but *not* the other way round: slowing down conflicts with direction). To construct them, it suffices to consider the particular case $\varphi = \mathrm{id}\uparrow\mathbf{I}$ (and compose with an arbitrary homotopy); thus

(4)
$$\begin{array}{ccc} f & \xrightarrow{\varphi} & g \\ 0_f \downarrow & \Theta'' & \downarrow 0_g \\ f & \xrightarrow[\varphi+0]{} & g \end{array} \qquad \begin{array}{l} \varphi = \mathrm{id}\uparrow\mathbf{I}, \quad f = \partial^-, \quad g = \partial^+: \{*\} \to \uparrow\mathbf{I}, \\ \varphi(t) = t, \quad (\varphi + 0_g)(t) = (2t)\wedge 1, \\ \Theta''(t, t') = (1-t').t + t'.((2t)\wedge 1). \end{array}$$

In fact, $\Theta''$ is a linear interpolation (in t') *from* $\varphi$ *to* $\varphi + 0_g$; since $\varphi(t) \leq (\varphi + 0_g)(t)$, $\Theta''$ preserves the order of the square and is a d-map $\uparrow\mathbf{I}^2 \to \uparrow\mathbf{I}^2$.

(c) *Folding*. A double homotopy $\Phi: A \times \uparrow\mathbf{I}^2 \to X$ with faces $\varphi, \psi, \sigma, \tau$ (as below) produces a 2-homotopy $\Psi$, by pasting $\Phi$ with two double homotopies of connection (denoted by #)

(5)
$$\begin{array}{ccccccc} f & = & f & \xrightarrow{\sigma} & h & \xrightarrow{\psi} & g \\ \| & \# & \varphi \downarrow & \Phi & \downarrow \psi & \# & \| \\ f & \xrightarrow[\varphi]{} & k & \xrightarrow[\tau]{} & g & = & g \end{array} \qquad \Psi: 0_f + \sigma + \psi \to \varphi + \tau + 0_g: f \to g$$

which, together with accelerations, shows that $\sigma + \psi \simeq_2 \varphi + \tau$ (2.5).

**2.7. Controlling deformation.** Directed homotopy equivalence and deformation retracts (2.4) can be controlled, *step by step*.

We shall speak of an (immediate) *d-homotopy equivalence in the future* when both composed maps can be reached, from the identities, in one deformation step, *at the end of it* (t = 1)

(1)   $f: X \rightleftarrows Y : g,$  $\qquad\qquad\qquad \varphi: \mathrm{id}X \to gf, \quad \psi: \mathrm{id}Y \to fg;$

of an (immediate) *d-homotopy equivalence in the past* in the reflection-dual case: both homotopies start from the composed maps; and of an *n-step homotopy equivalence* for a sequence of n immediate equivalences (of any type)

(2)   $X = X_0 \rightleftarrows X_1 \rightleftarrows \ldots \rightleftarrows X_n = Y;$

composing them, it is easy to see that $X \simeq Y$, as already defined (2.4).

Similarly, a subspace $X_0$ of a d-space $X$ will be said to be an (immediate) *future* (resp. *past*) *deformation retract* of X if its inclusion u has a retraction p with $\mathrm{id}X \preceq up$ (resp. $\mathrm{id}X \succeq up$). We say that $X_0$ is a *deformation retract in* n *steps* if there is a finite sequence

(3)   $X_0 \subset X_1 \subset \ldots \subset X_n = X,$

where each d-space is an immediate deformation retract of the following one; and in *precisely n steps* if a shorter chain does not exist. Again, composing all retractions, it is easy to see that $X_0$ is a deformation retract of $X$ as previously defined (2.4), and that $X_0$ and $X$ are d-homotopy equivalent in $n$ steps. Note also that, in (3), each $X_i$ is either *upper bounded* by $X_{i-1}$ (each point of $X_i$ has some upper bound in $X_{i-1}$) or *lower bounded* there, for the path preorder.

Plainly, $X_0 \subset X$ is a future deformation retract of $X$ iff there is a d-map $\varphi$ such that

(4) $\quad \varphi: X \times {\uparrow}\mathbf{I} \to X, \qquad\qquad \varphi(x, 0) = x, \quad \varphi(x, 1) \in X_0 \qquad\qquad (x \in X);$

for instance, the cylinder $X \times {\uparrow}\mathbf{I}$ has a strong future deformation retract at its upper basis $\partial^+(X) \subset {\uparrow}\mathbf{I}X$, by the lower connection (reaching the upper basis at time $t' = 1$)

(5) $\quad g^-: (X \times {\uparrow}\mathbf{I}) \times {\uparrow}\mathbf{I} \to X \times {\uparrow}\mathbf{I}, \qquad\qquad g^-(x, t, t') = (x, t \vee t'),$

and a strong past deformation retract at its lower basis $\partial^-(X)$.

A d-space $X$ is *future contractible*, or *past contractible*, or *contractible in* $n$ *steps*, if it has such a deformation retract, at a point. $X$ is future contractible iff there is a homotopy $\mathrm{id}X \to v$ with a constant map $v: X \to X$; then $x^+ = v(X)$ is *a* maximum for the path preorder $\preceq$.

The cones ${\uparrow}C^+(S^1)$, ${\uparrow}C^+({\uparrow}S^1)$ (1.6.2-3) are future contractible to their vertex, while the lower ones are past contractible. The half-line ${\uparrow}]-\infty, 0] \subset {\uparrow}\mathbf{R}$ is future contractible to $0$ (which is reached by the d-homotopy $\varphi(x, t) = (1-t)x$, at time $t = 1$), but not past contractible. The d-line ${\uparrow}\mathbf{R}$ is 2-step contractible, as ${\uparrow}]-\infty, 0]$ is a past deformation retract (reached by the homotopy $\psi(x, t) = x \wedge (tx)$ at $t = 0$); two steps are needed, since the line has no extreme for the (path) order.

Similarly, all ${\uparrow}\mathbf{R}^n$ are precisely contractible in 2 steps $(n > 0)$: one can take as a past deformation retract $X = {\uparrow}]-\infty, 0]^n$, moving all points of the complement to the boundary of $X$ along lines parallel to the main diagonal $x_1 = ... = x_n$. The V-shaped d-space $V$

(6) $\quad V = ([0, 1] \times \{0\}) \cup (\{0\} \times [0, 1]) \subset {\uparrow}\mathbf{R}^2,$

is past contractible (to the origin). The *infinite stairway* $W$

(7) $\quad W = \cup_{k \in \mathbf{Z}} (([k, k+1] \times \{-k\}) \cup (\{k\} \times [-k, 1-k])) \subset {\uparrow}\mathbf{R}^2,$

is not d-contractible (in fact, $f \preceq \mathrm{id}W$ or $\mathrm{id}W \preceq f$ imply $f(W) = W$). A *finite stairway* consisting of $2n$ or $2n-1$ consecutive segments of $W$ is contractible in $n$ steps (in the even case, each step contracts the first and last segment; in the odd one, the first step contracts one of them).

## 3. Computing the fundamental category

The fundamental category of a d-space is introduced. Non obvious computations are based on a van Kampen-type theorem (3.6), similar to R. Brown's version for the fundamental groupoid of spaces [Br].

**3.1. The fundamental category.** Directed paths are now considered modulo 2-homotopy, i.e. homotopy with fixed endpoints.

A *double path* in $X$ is a d-map $A: {\uparrow}\mathbf{I}^2 \to X$. It is the elementary instance of a double homotopy (2.5), defined on the point, and the previous results apply; its four faces are paths in $X$, between four vertices. A *2-path* is a double path whose faces $\partial_1^\alpha$ are degenerate; it is a 2-homotopy A: a $\preceq_2$ b: x $\to$ x' between its faces $\partial_2^\alpha$, which have the same endpoints. A 2-homotopy class of paths [a] is a class of the equivalence relation $\simeq_2$ spanned by the preorder $\preceq_2$.

The *fundamental category* $\uparrow\Pi_1(X)$ of a d-space has for objects the points of X; for arrows [a]: x $\to$ x' the 2-homotopy classes of paths from x to x', as defined above. Composition – written additively – is induced by concatenation of consecutive paths, identities derive from degenerate paths

(1)    [a] + [b]  =  [a+b],               $0_x$ = [e(x)] = $[0_x]$.

We prove below that $\uparrow\Pi_1(X)$ is indeed a category, and that the obvious action on arrows defines a functor $\uparrow\Pi_1$: d**Top** $\to$ **Cat**  (small categories)

(2)    $\uparrow\Pi_1(f)(x)$ = f(x),               $\uparrow\Pi_1(f)[a]$ = $f_{*1}[a]$ = [fa],

invariant with respect to notions of directed homotopy in **Cat** which will be developed in the next Section (4.1). The fundamental category of X is linked to the fundamental groupoid of the underlying space UX, by the obvious *comparison* functor

(3)    $\uparrow\Pi_1(X) \to \Pi_1(UX)$,           x $\mapsto$ x,    [a] $\mapsto$ [a],

which is the identity on objects, but need not be full (obviously) nor faithful (3.5). Plainly, if X is a topological space with the natural d-structure (X = $C^0$UX), then $\uparrow\Pi_1(X) = \Pi_1(UX)$.

**3.2. Invariance Theorem.** (a) For every d-space X, $\uparrow\Pi_1(X)$ is a category and the previous formulae (3.1.2) do define a functor, which preserves sums and products. The reflected d-space gives the opposite category, $\uparrow\Pi_1(RX) = (\uparrow\Pi_1(X))^{op}$.

(b) If a: x $\to$ x' is a reversible path (1.5), its class [a] is an invertible arrow in $\uparrow\Pi_1(X)$.

(c) The functor $\uparrow\Pi_1$: d**Top** $\to$ **Cat** is d-homotopy invariant, in the following sense: a d-homotopy $\varphi$: f $\to$ g: X $\to$ Y induces a natural transformation (a d-homotopy of categories, 4.1)

(1)    $\varphi_{*1}$: $f_{*1} \to g_{*1}$: $\uparrow\Pi_1(X) \to \uparrow\Pi_1(Y)$,        $\varphi_{*1}(x)$ = [$\varphi(x)$]: f(x) $\to$ g(x),

where $\varphi(x)$ is the path in Y, from the representative map X $\to$ $\uparrow$PY. (Therefore $\uparrow\Pi_1$ *preserves d-homotopy*, *d-homotopy equivalence and deformation retracts*, in the sense of 4.1).

(d) A *reversible* d-homotopy $\varphi$ induces an *invertible* transformation $\varphi_{*1}$. (Therefore $\uparrow\Pi_1$ turns *reversible* d-homotopy equivalence into equivalence of categories).

**Proof.** (a) Composition is well defined, in 3.1.1. Given 2-homotopies A: a $\preceq_2$ a': x $\to$ x' and B: b $\preceq_2$ b': x' $\to$ x", the pasting A $+_1$ B: a+b $\preceq_2$ a'+b': x $\to$ x" shows that [a+b] = [a'+b']. The general case, for the equivalence relation $\simeq_2$, follows by taking, in A or B, a trivial 2-homotopy and applying transitivity. The fact that a d-map f: X $\to$ Y gives a well-defined transformation $\uparrow\Pi_1(f)[a]$ = [fa] is also obvious: for A: a $\preceq_2$ a', take fA: fa $\preceq_2$ fa'.

In $\uparrow\Pi_1(X)$, constant paths produce (strict) identities, because of the acceleration 2-homotopies $0_x$+a $\to$ a $\to$ a+$0_{x'}$ (2.6.3). Associativity, on three consecutive paths a, b, c in X, follows from considering a 2-homotopy B: (0+a)+(b+c) $\to$ (a+b)+(c+0), constructed by pasting double paths deriving from degeneracy and connections (all denoted by #)

(1)
```
x ═══ x ──a→ y ──b→ z ──c→ w
‖   #  ↓a #  ‖   #  ‖   #  ‖
x ──→ y ═══ y ──→ z ──→ w
‖   #  ‖   #  ↓b #  ‖   #  ‖
x ──→ y ──→ z ═══ z ──→ w
‖   #  ‖   #  ‖   #  ↓c #  ‖
x ──a→ y ──b→ z ──c→ w ═══ w
```

The argument is concluded by two other 2-homotopies, deriving from accelerations; they cannot be pasted with B, because of conflicting directions

(2)    A: $(a+b)+c \to (a+b)+(c+0)$,          C: $(0+a)+(b+c) \to a+(b+c)$.

The preservation of sums and cartesian products by the functor $\uparrow\Pi_1$ is proved in the same (easy) way as in the ordinary case.

(b) By definition, $a: \mathbf{I}^\sim \to X$ is assumed to be a d-map. The double path

(3)    $A = ag^\sim.(\mathbf{I}^\sim \times r): \mathbf{I}^{\sim 2} \to X$,

$$\begin{array}{ccc} x' & \xrightarrow{0} & x' \\ a^{op} \downarrow & A & \downarrow 0 \\ x & \xrightarrow{a} & x' \end{array}$$

is indeed directed with respect to the reversible structures; in fact, given two piecewise monotone real functions h, k, also h∨k is so (if h is increasing and k decreasing on some interval $[t_0, t_1]$, and $h(t) = k(t)$ at some intermediate point, then h∨k coincides with k on $[t_0, t]$, with h on $[t, t_1]$). Finally, by folding (2.6c), and recalling that $\uparrow\mathbf{I}$ is finer than $\mathbf{I}^\sim$, we get a 2-path

(4)    A': $\uparrow\mathbf{I}^2 \to \mathbf{I}^{\sim 2} \to \mathbf{I}^\sim \to X$,          A': $0 \to a^{op} + a + 0: \uparrow\mathbf{I} \to X$.

(c) The naturality of the transformation associated to $\varphi: f \to g$ on the arrow $[a]: x \to x'$ in $\uparrow\Pi_1(X)$ amounts to the relation $[fa] + [\varphi(x')] = [\varphi(x)] + [ga]$. This follows from the existence of the double path $\Phi = \varphi \circ a = \varphi.(a \times \uparrow\mathbf{I}): \{*\} \times \uparrow\mathbf{I}^2 \to X$

(5)
$$\begin{array}{ccc} f(x) & \xrightarrow{fa} & f(x') \\ \varphi(x) \downarrow & \varphi \circ a & \downarrow \varphi(x') \\ g(x) & \xrightarrow{ga} & g(x') \end{array} \qquad fa + \varphi(x') \preceq_2 \varphi(x) + ga.$$

(d) Is a straightforward consequence of (b) and (c).

**3.3. Homotopy monoids.** The *fundamental monoid* $\uparrow\pi_1(X, x)$ of the d-space X at the point x is the monoid of endoarrows $[c]: x \to x$ in $\uparrow\Pi_1(X)$. It forms a functor from the (obvious) category d**Top**$_*$ of *pointed d-spaces*, to the category of monoids

(1)    $\uparrow\pi_1$: d**Top**$_* \to$ **Mon**,          $\uparrow\pi_1(X, x) = \uparrow\Pi_1(X)(x, x)$,

which is *strictly* d-homotopy invariant: a *pointed d-homotopy* $\varphi: f \to g: (X, x) \to (Y, y)$ has, by definition, a trivial path at the base-point ($\varphi(x) = 0_y$), whence $f_{*1} = g_{*1}$ (3.2.5).

Similarly, we have a functor from the comma category $d\mathbf{Top}\backslash \mathbf{S}^0$ of *bipointed d-spaces*

(2)   $\uparrow \pi_1: d\mathbf{Top}\backslash \mathbf{S}^0 \to \mathbf{Set}$,     $\uparrow \pi_1(X, x, x') = \uparrow \Pi_1(X)(x, x')$,

which is strictly d-homotopy invariant, up to *bipointed d-homotopies* (leaving fixed each base point). One can view (1) and (2) as representable homotopy functors (2.4) on $d\mathbf{Top}_*$ and $d\mathbf{Top}\backslash \mathbf{S}^0$, which 'accounts' for their strict invariance. Moreover, both can be computed by the methods developed below for $\uparrow \Pi_1 X$. (For the homotopy structure of comma categories, see [G3].)

The existence of a *reversible* path (1.5) from $x$ to $x'$ implies that their fundamental monoids are isomorphic (by 3.2b); without reversibility, this need not be true (cf. 3.5). However, in a *homogeneous* d-space, where $\mathrm{Aut}(X)$ acts transitively, all $\uparrow \pi_1(X, x)$ are plainly isomorphic; this applies, for instance, to the directed circle $\uparrow \mathbf{S}^1$.

**3.4. Simple d-spaces.** Say that a d-space $X$ is *1-simple* if its fundamental category is a preorder, or equivalently if $\uparrow \Pi_1 X = \mathrm{cat}(X, \preceq)$, the category associated with the path preorder:

(1)   $\uparrow \Pi_1 X(x, x')$ has one arrow when $x \preceq x'$, no arrow otherwise.

(a) Every *convex* subset $X$ of $\mathbf{R}^n$, with the order structure induced by $\uparrow \mathbf{R}^n$, is 1-simple. In fact, if $x \leq x'$, there is a d-path from $x$ to $x'$, e.g. $a(t) = (1-t).x + t.x'$; the converse is obvious. Moreover, given two increasing paths $a, b: \uparrow \mathbf{I} \to X$ from $x$ to $x'$, we can always assume that $a \leq b$ (otherwise, replace the first with $0_x + a \preceq_2 a$, the second with $b + 0_{x'} \succeq_2 b$); then, the interpolation 2-path $A(t, t') = (1-t').a(t) + t'.b(t)$ preserves the order of $\uparrow \mathbf{I}^2$ and provides a directed 2-homotopy $A: a \to b$.

(b) It follows that a d-space $X$ is certainly 1-simple whenever the following condition holds: if $x' \preceq x''$, then the d-subspace $\{x \in X \mid x' \preceq x \preceq x''\}$ is isomorphic to some convex d-subspace of $\uparrow \mathbf{R}^n$.

(c) The following objects are 1-simple: any non-empty interval $J \subset \uparrow \mathbf{R}$; any product of such in $\uparrow \mathbf{R}^n$; $V, W \subset \uparrow \mathbf{R}^2$ (2.7.6-7); any 'fan' formed by the union of (finitely or infinitely many) segments or half-lines spreading from a point, in some $\uparrow \mathbf{R}^n$. (Here, one should not confuse the path-order with the order induced by $\uparrow \mathbf{R}^n$, which is coarser and of less interest.)

(d) The ordered circle $\uparrow \mathbf{O}^1$ is not 1-simple. Any d-path there stays either in the left half or in the right one, whence the two obvious d-paths moving in such half-circles $a, b: x^- \to x^+$ (from the minimum $x^- = (0, -1)$ to the maximum $x^+ = (0, 1)$) are not 2-homotopic. Adding this consideration to the previous ones, it is easy to determine the fundamental category

(2)   $\uparrow \Pi_1 \uparrow \mathbf{O}^1(x, x') = \{[a], [b]\}$         if $x = x^-$, $x' = x^+$         ($[a] \neq [b]$),

and, otherwise, one arrow if $x \preceq x'$, no arrows in the contrary.

**3.5. Comments.** The fundamental category contains information which can disappear modulo d-homotopy, yet *not* modulo (bi)pointed d-homotopy (3.3) nor *reversible* d-homotopy (3.2d). This can be traced back to the fact that, for a d-homotopy equivalence $p: X \to Y$, the induced functor $p_{*1}: \uparrow \Pi_1 X \to \uparrow \Pi_1 Y$ need not be *full* nor *faithful* (even when $p$ is a strong deformation retraction).

For the first case, just consider the fact that ↑**I** is strongly past contractible to 0, while ↑$\Pi_1$(↑**I**) = cat(**I**, ≤) keeps the information of the (path) order; thus $p_{*1}$: ↑$\Pi_1$(↑**I**) → ↑$\Pi_1${0} is not full.

The infraction to faithfulness can produce even more unusual effects (but we shall see that it cannot happen if ↑$\Pi_1 X$ satisfies the cancellation laws, 4.3a):

(i) a strongly d-contractible object X can have loops c which are not homotopically trivial ([c] ≠ 0),

(ii) such loops are then annihilated by the deformation retraction p: X → {∗} ($p_{*1}$ is not faithful),

(iii) such loops are 'loop-homotopic' to the constant loop, without being 2-homotopic to it.

In fact, take the disc with the structure X = ↑$C^-$(**S**$^1$) (1.6.2), which is strongly past contractible to its centre $v^-$, but not reversibly so. Any concentric circle C inherits the natural structure of **S**$^1$, and no path between two of its points x', x" can leave it; thus, the restriction of ↑$\Pi_1 X$ to the points of C coincides with the fundamental groupoid of the circle and has d-loops c: **S**$^1$ → X with [c] ≠ 0. Any deformation φ: X×↑**I** → X with φ(x, 0) = $v^-$, φ(x, 1) = x (2.7.4) yields a loop-homotopy φ.(c×↑**I**): **S**$^1$×↑**I** → X from $0_{v^-}$ to c. Note also that, if c is a loop at $x_0$, the homotopy class of the path a = φ.($x_0$, –): $v^-$ → $x_0$ is not cancellable in ↑$\Pi_1 X$ (not epi): [a] + [c] = [a].

Similar arguments also allow us to completely determine the fundamental category of ↑$C^-$(**S**$^1$): from the origin to any point there is *one* arrow (this would follow directly from 4.1: $v^-$ must be initial in ↑$\Pi_1 X$); otherwise, there are arrows determined by their 'winding number' around the origin

(1)    ↑$\Pi_1(X)(x, x')$ = $\Pi_1 \mathbf{S}^1(q(x), q(x'))$                          (0 < ||x|| ≤ ||x'||),

where q(x) = x/||x||; there are no other arrows; concatenation works accordingly, as in $\Pi_1 \mathbf{S}^1$.

**3.6. Pasting Theorem** ('Seifert-van Kampen' for fundamental categories). Let X be a d-space; $X_1, X_2$ two d-subspaces and $X_0 = X_1 \cap X_2$.

(a) If X = int($X_1$)∪int($X_2$), the following diagram of categories and functors (induced by inclusions) is a pushout in **Cat**

(1)
$$\begin{array}{ccc} \uparrow\Pi_1 X_0 & \xrightarrow{u_1} & \uparrow\Pi_1 X_1 \\ u_2 \downarrow & {}_{v_2} & \downarrow v_1 \\ \uparrow\Pi_1 X_2 & \longrightarrow & \uparrow\Pi_1 X \end{array}$$

(b) More generally, the same fact holds provided one can find two d-subspaces $w_i$: $Y_i \subset X_i$ with retractions $p_i$: $Y_i$ → $X_i$ (d-maps with $p_i w_i = idY_i$; no deformation is required) such that:

(2)    X = int($Y_1$)∪int($Y_2$),                           $p_1$ and $p_2$ coincide on $Y_0 = Y_1 \cap Y_2$.

**Proof.** (a) We shall use the n-ary concatenation of consecutive d-paths, written $a_1 + ... + a_n$ (1.1). Let $F_i$: ↑$\Pi_1 X_i$ → C be two functors which coincide on ↑$\Pi_1 X_0$ ($F_1 u_1 = F_2 u_2$); we have to prove that they have a unique 'extension' F: ↑$\Pi_1 X$ → C. On the objects, this is obvious since |X| = |$X_1$| ∪ |$X_2$| and |$X_0$| = |$X_1$| ∩ |$X_2$|.

Let then a: x → y: {∗} → X be a path. By the Lebesgue covering lemma, there is a finite decomposition 0 < 1/n < 2/n ... < 1 of the standard interval such that each subinterval [(i–1)/n, i/n] is mapped by a into $X_1$ or $X_2$ (a *suitable* decomposition for our data). Thus, a = $a_1$ + ... + $a_n$ where each $a_i: [0, 1] \to X$ is a directed path (by increasing reparametrisation) contained in some $X_{k_i}$, hence a d-path there. Define (using the additive notation also for composition in C)

(3) $\quad F[a] = F_{k_1}[a_1] + ... + F_{k_n}[a_n] \in C(F(x), F(x'))$.

First, this does not depend on choosing $k_i$: if $Im(a_i) \subset X_1 \cap X_2 = X_0$, then $F_1 u_1 = F_2 u_2$ shows that $F_1[a_i] = F_2[a_i]$. Second, this does not depend on the choice of n: if also m gives a suitable partition, use the partition deriving from mn to prove that they give the same result. Third, F[a] does not depend on the representative path a. It is sufficient to show this for a second path a': x → x', linked to the first by a 2-path A: a → a'; in other words, A: $\uparrow\mathbf{I}^2 \to X$ has degenerate 1-faces, and 2-faces coinciding with a, a'. Again by the Lebesgue covering lemma, applied to the compact metric square $[0, 1]^2$, there is some integer n > 0 such that all elementary squares $[(i–1)/n, i/n] \times [(j–1)/n, j/n]$ are mapped by A into $X_1$ or $X_2$. A can be obtained as an 'n×n-pasting' of its (reparametrised) restrictions to these squares, $A_{ij}: \uparrow\mathbf{I}^2 \to X_{k(i,j)} \subset X$

(4) $\quad A = (A_{11} +_1 A_{21} +_1 ... +_1 A_{n1}) +_2 ... +_2 (A_{1n} +_1 A_{2n} +_1 ... +_1 A_{nn})$.

Every square $B = A_{ij}$ produces, by folding (2.6c), a 2-homotopy relation in $X_{k(i,j)}$

(5) $\quad \partial_2^- B + \partial_1^+ B \simeq_2 \partial_1^- B + \partial_2^+ B$.

Therefore, using the fact that all 1-faces on the boundary are degenerate ($\partial_1^- A_{1i}, \partial_1^+ A_{ni}$), and the coincidence of faces between contiguous little squares, we can gradually move from a to a'

(6) $\quad F[a] = F_{k(1,1)}[\partial_2^- A_{11}] + ... + F_{k(n,1)}[\partial_2^- A_{n1}]$

$\quad\quad = F_{k(1,1)}[\partial_2^- A_{11}] + ... + F_{k(n,1)}[\partial_2^- A_{n1}] + F_{k(n,1)}[\partial_1^+ A_{n1}]$ (by degeneracy)

$\quad\quad = F_{k(1,1)}[\partial_2^- A_{11}] + ... + F_{k(n,1)}[\partial_1^- A_{n1}] + F_{k(n,1)}[\partial_2^+ A_{n1}]$ (by (5))

$\quad\quad = F_{k(1,1)}[\partial_2^- A_{11}] + ... + F_{k(n,1)}[\partial_1^+ A_{n-1,1}] + F_{k(n,2)}[\partial_2^- A_{n2}]$ (by contiguity)

$\quad\quad = ... = F_{k(1,2)}[\partial_2^- A_{21}] + ... + F_{k(n,2)}[\partial_2^- A_{n1}]$

$\quad\quad = ... = F_{k(1,n)}[\partial_n^- A_{n1}] + ... + F_{k(n,n)}[\partial_n^- A_{n1}] = F[a']$.

Thus, F: $\uparrow\Pi_1 X \to C$ is also well defined on arrows. To show that it preserves composition just note that, if two consecutive d-paths a, b have a suitable decomposition on n subintervals, than a+b inherits a suitable decomposition $a + b = a_1 + ... + a_n + b_1 + ... + b_n$ which keeps the original paths apart. Finally, the uniqueness of the functor F is obvious.

(b) By (a), the square deriving from $Y_0, Y_1, Y_2$, and $Y = X$ is a pushout of categories.

Also the inclusion $w_0: X_0 \subset Y_0$ has a retraction $p_0$, the common restriction of $p_1$ and $p_2$ to $Y_0$. Therefore, all $w_i$ and $p_i$ form a retraction in the category of commutative squares of d**Top**

(7) $\quad \mathbf{w} = (w_0, w_1, w_2, idX): \mathbf{X} \to \mathbf{Y}, \quad\quad \mathbf{p} = (p_0, p_1, p_2, idX): \mathbf{Y} \to \mathbf{X} \quad (\mathbf{pw} = idX)$,

The functor $\uparrow\Pi_1$ takes all this into a retraction $\mathbf{w}_*: \uparrow\Pi_1\mathbf{X} \rightleftarrows \uparrow\Pi_1\mathbf{Y} :\mathbf{p}_*$ in the category of commutative squares of **Cat**. Since $\uparrow\Pi_1\mathbf{Y}$ is a pushout, also its retract $\uparrow\Pi_1\mathbf{X}$ is so (as can be easily checked, or seen in [Br], 6.6.7).

**3.7. Computations.** Putting together the pasting theorem 3.6 and the preceding results, it is easy to compute the fundamental category of the directed circle, of the d-spaces of the Introduction, etc.

(a) The directed circle ↑**S**¹. Apply van Kampen (3.6a) in the obvious way, with two arcs $X_1$, $X_2$ isomorphic to ↑**I** and $X_0 \cong$ ↑**I** + ↑**I**. The resulting pushout in **Cat** shows that ↑$\Pi_1$↑**S**¹ is the subcategory of the groupoid $\Pi_1$**S**¹ formed by the classes of anticlockwise paths. In particular, each monoid ↑$\pi_1$(↑**S**¹, x) is isomorphic to the additive monoid **N** of natural numbers.

(b) The fundamental category of ↑C⁻(↑**S**¹) can now be easily deduced, as for ↑C⁻(**S**¹) in 3.5.

(c) The (compact) d-space H ⊂ ↑**R**² represented below

(1) 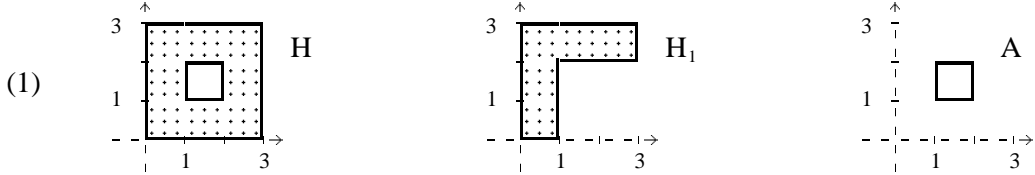

has a fundamental category 'similar' to that of the ordered circle (3.4d)

(2)    ↑$\Pi_1$H(x, x'):        two arrows if x ∈ [0, 1]² and x' ∈ [2, 3]²,

and otherwise, one arrow if $x \preceq x'$, no arrows in the contrary. This follows from applying van Kampen with $H_1$ as in (1) and a similar $H_2$; they are 1-simple spaces, again by 3.6.

Note that the d-subspace A is a strong deformation retract of H (in 2 steps), isomorphic to the ordered circle. However, while this determines ↑$\Pi_1$A as the full subcategory of ↑$\Pi_1$H with objects in A (by 4.3b), it is not clear how to deduce ↑$\Pi_1$H from ↑$\Pi_1$A.

(d) One can similarly compute the fundamental category of d-spaces like the ones of the Introduction, fig. (2). Again, there are facts which can disappear up to d-homotopy: for instance, at the left of each obstruction in X or Y there is a triangle from which no path can reach b; which is of interest in all the interpretations we were considering. This information, however, is stable under *bipointed* d-homotopy (3.3.2), with one point in that triangle and the second at b.

## 4. Complements

Directed homotopy of categories is briefly dealt with and related to d-spaces. We end by discussing directed geometric realisation of cubical sets (cf. [BrH]) and directed metrisability of d-spaces.

**4.1. Directed homotopy for categories.** Classically, the homotopical invariance of the fundamental groupoid functor $\Pi_1$: **Top** → **Gpd** (small groupoids) means that it turns homotopy equivalence of spaces into ordinary equivalence of groupoids; the latter can be viewed as homotopy equivalence in **Gpd**, based on a reversible *interval*, the groupoid **i** = {0 ⇌ 1}. Now, for ↑$\Pi_1$, we need to replace equivalence of groupoids by a directed notion.

We shall give a brief description of *directed homotopy in* **Cat** (the category of small categories). This will be based on the *directed interval* ↑**i** = **2** = {0 → 1}, an order category, with obvious faces $\partial^{\pm}$: **1** → **2**, where **1** = {0} is the pointlike category. Again, **2** is an internal lattice (2.1), and reversion is missing but partially surrogated by the reflecting isomorphism, r: **2** → **2**ᵒᵖ.

A *point* x: **1** → C is an object of C. A *d-path* a: **2** → C from x to x' is an arrow a: x → x' of C; their concatenation is defined by the composition in C, and is strictly associative, with strict identities. To be consistent with the previous notation, we shall write additively, a + b, the composition of consecutive arrows (in the objects of **Cat**) and $0_x$ the identity at x. (Formally, we should note that the concatenation pushout 1.5.1 gives here the order category **3**, and not the interval **2**; concatenation is realised by an obvious functor k: **2** → **3**.) A double path **2**×**2** → C is a commutative square, while a 2-path is necessarily trivial. Therefore $\uparrow\Pi_1 C$, defined as above for d-spaces, just *coincides* with C, and $\uparrow\Pi_1 f = f$ on functors. A *reversible path* a: **i** → C is an isomorphism.

The directed *cylinder* $\uparrow IC = C \times \mathbf{2}$ and its right adjoint, $\uparrow PD = D^{\mathbf{2}}$ (the category of morphisms of D) show that a *d-homotopy* φ: f → g: C → D is the same as a natural transformation between functors. Operations between d-homotopies and functors are defined as previously and amount to the usual operations; concatenation is vertical composition, written as (φ+ψ)(x) = φx + ψx; it is strictly associative, with strict identities $0_f$: f → f, the vertical identities of functors.

Given two parallel functors f, g: C → D (and proceeding as in 2.3), we write f $\preceq$ g if there is a natural transformation f → g, and f ≃ g if there is a finite sequence of them, f = $f_0$ → $f_1$ ← $f_2$ ... $f_n$ = g. The categories C and D are *d-homotopy equivalent* (C ≃$_d$ D) if there are functors

(1)   f : C ⇄ D : g,                              gf ≃ idC,    fg ≃ idD;

C is a *deformation retract* of D if, moreover, gf = idC. Finally, C is *d-contractible* if it is d-homotopy equivalent to **1**. All this can be controlled step-by-step, as in 2.7. It is easy to see that a category has an initial object (say v⁻) iff it is strongly past contractible (to v⁻), iff the functor p: C → **1** has a left adjoint (v⁻: **1** → C); thus, all ordinals > 0 are d-homotopy equivalent.

Directed homotopy equivalence distinguishes new 'shapes', in **Cat**. It is weaker than ordinary equivalence, which corresponds to *reversible homotopies*, based on **i**. But it implies ordinary homotopy equivalence of the classifying spaces, the geometric realisations of the simplicial nerves (which is an undirected notion): a natural transformation φ: f → g: C → D produces, by the nerve functor N: **Cat** → **Smp**, a simplicial homotopy Nφ: Nf → Ng: NC → ND (because N(C×**2**) = NC×N**2** and N**2** is the simplicial interval •—→•), and then, by ordinary geometric realisation, a homotopy of the classifying spaces (cf. [May], § 16; [Cu] 1.29).

**4.2. Lemma.** Let φ: h → k: C → D be a natural transformation all whose components φx (x ∈ ObC) are cancellable in D (mono and epi). Then h is faithful if and only if k is so.

**Proof.** Take two arrows $u_1$, $u_2$: x → x' in C and the resulting commutative squares in D

(1)   
$$\begin{array}{ccc} h(x) & \xrightarrow{\varphi x} & h(x) \\ {}_{hu_i}\downarrow & & \downarrow {}^{ku_i} \\ k(x') & \xrightarrow{\varphi x'} & h(x') \end{array}$$

If h is faithful and k($u_1$) = k($u_2$), cancelling φx' (mono) we obtain h($u_1$) = h($u_2$) and $u_1$ = $u_2$. Conversely, if k is faithful, use the fact that φx is epi.

**4.3. Theorem** (Distinguishing d-homotopy). (a) Let $f: X \to Y$ be a d-homotopy equivalence in d**Top** (or in **Cat**). If $\uparrow\Pi_1 X$ satisfies the cancellation laws (all arrows are mono and epi), then $f_{*1}$: $\uparrow\Pi_1 X \to \uparrow\Pi_1 Y$ is faithful.

(b) Let $u: X \subset Y$ be a *strong* deformation retract in d**Top** (or in **Cat**). Then $u_{*1}: \uparrow\Pi_1 X \to \uparrow\Pi_1 Y$ is a full embedding, while its retraction $p_{*1}$ need not be faithful nor full.

(c) Let $h \simeq k: X \to Y$ in d**Top**. If all the paths $\varphi_i(x): h_{i-1}(x) \to h_i(x)$ which intervene in the previous relation are cancellable *in* $\uparrow\Pi_1 Y$ (i.e., the corresponding classes $[\varphi_i(x)]$ are mono and epi), then the functor $h_{*1}: \uparrow\Pi_1 X \to \uparrow\Pi_1 Y$ is faithful if and only if $k_{*1}$ is so.

**Proof.** First, (c) follows from the previous lemma. Recall now that, in **Cat**, $\uparrow\Pi_1 X = X$ (4.1).

(a) Take a d-map $g: Y \to X$ with $gf \simeq idX$. By (c), $g_{*1}f_{*1}$ is faithful, whence also $f_{*1}$ is so.

(b) We already know that $p_{*1}$ need not be faithful nor full (3.5). By hypothesis, all the d-homotopies which intervene in the relation $up \simeq id: Y \to Y$ can be chosen to be trivial on $X$ (2.4). By (c), the functor $u_{*1}p_{*1}: \uparrow\Pi_1 Y \to \uparrow\Pi_1 Y$ is faithful *on the points of* $X$ (i.e., on the full subcategory of such objects in $\uparrow\Pi_1 Y$), and also $p_{*1}: \uparrow\Pi_1 Y \to \uparrow\Pi_1 X$ is so. Now, given two points $x, x' \in X$, consider the retraction

(1)     $u_{*1}: \uparrow\Pi_1 X(x, x') \rightleftarrows \uparrow\Pi_1 Y(x, x') : p_{*1}$     $(p_{*1}.u_{*1} = id)$;

since we have already proved that this mapping $p_{*1}$ is injective, $u_{*1}$ must be its (bilateral) inverse.

**4.4. Applications to d-spaces.** We have essentially proposed d-homotopy to study directed paths. The fundamental category can also be used to distinguish the d-homotopy type of d-spaces (when the underlying spaces are homotopy equivalent).

Thus, to show that $S^1$, $\uparrow S^1$ and $\uparrow O^1$ are not d-homotopy equivalent, one can apply 4.3a (after computing their fundamental categories, in 3.4 and 3.7): there are no faithful functors from $\Pi_1 S^1$ to $\uparrow\Pi_1\uparrow S^1$ or $\uparrow\Pi_1\uparrow O^1$, nor from $\uparrow\Pi_1\uparrow S^1$ to $\uparrow\Pi_1\uparrow O^1$; and the cancellation laws hold. (But we have already proved the same in 2.4, resting on ordinary homotopy and elementary arguments on d-maps).

It is easy to realise, on the 'eight figure' $X$, locally preordered d-structures with $\uparrow\pi_1(X, x)$ isomorphic to $\mathbf{Z}*\mathbf{Z}$, $\mathbf{Z}*\mathbf{N}$, $\mathbf{N}*\mathbf{N}$ ($*$ is the 'free direct product' of monoids, i.e. their categorical sum)

(1) 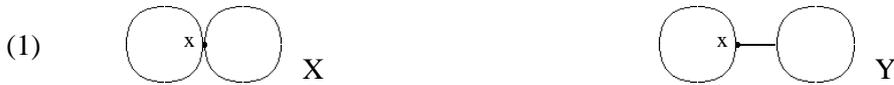

On $Y$ one can obtain the same, and also $\mathbf{Z}\times\mathbf{Z}$, $\mathbf{Z}\times\mathbf{N}$, $\mathbf{N}\times\mathbf{N}$ (ordering the segment); all these structures are distinguished by $\uparrow\Pi_1$.

**4.5. Directed geometric realisation.** Cubical sets have a clear realisation as directed spaces, since we obviously want to realise the object with one free generator of degree n as $\uparrow\mathbf{I}^n$. On the other hand, for simplicial sets, it is not obvious which d-structure we should assign to the standard simplex $\Delta^n$ (and even less if we want this to be consistent with barycentric subdivision).

A cubical set (with faces and degeneracies) $K = ((K_n), (\partial_i^\alpha), (e_i))$ can be viewed as a functor $K: \mathbb{I}^{op} \to \mathbf{Set}$, on a category $\mathbb{I} \subset \mathbf{Set}$ (the *cubical site*). Its objects are the sets $\mathbf{2}^n = \{0, 1\}^n$, its mappings are generated by the elementary faces $\partial^\alpha: \mathbf{2}^0 \to \mathbf{2}$ and degeneracy $e: \mathbf{2} \to \mathbf{2}^0$, under finite products (in **Set**) and composition. Equivalently, the mappings of $\mathbb{I}$ are generated (under composition) by the following faces and degeneracies ($i = 1,..., n$; $\alpha = 0, 1$; $t_i = 0, 1$)

(1) $\quad \partial_i^\alpha = \mathbf{2}^{i-1}\times\partial^\alpha\times\mathbf{2}^{n-i}\colon \mathbf{2}^{n-1} \to \mathbf{2}^n, \qquad \partial_i(t_1,..., t_n) = (t_1,..., t_{i-1}, \alpha,..., t_n),$

$\quad\ e_i = \mathbf{2}^{i-1}\times e\times\mathbf{2}^{n-i}\colon \mathbf{2}^n \to \mathbf{2}^{n-1}, \qquad e_i(t_1,..., t_n) = (t_1,..., \hat{t}_i,..., t_n).$

There is an obvious embedding of $\mathbb{I}$ in d**Top**, where faces and degeneracies are realised as above, with $\partial^\alpha\colon \{*\} \rightrightarrows \uparrow\mathbf{I} \colon e$ (and $t_i \in \uparrow\mathbf{I}$)

(2) $\quad \uparrow I\colon \mathbb{I} \to \mathbf{dTop}, \qquad\qquad\qquad\qquad\quad \mathbf{2}^n \mapsto \uparrow\mathbf{I}^n,$

$\quad\ \partial_i^\alpha = \uparrow\mathbf{I}^{i-1}\times\partial^\alpha\times\uparrow\mathbf{I}^{n-i}\colon \uparrow\mathbf{I}^{n-1} \to \uparrow\mathbf{I}^n, \qquad e_i = \uparrow\mathbf{I}^{i-1}\times e\times\uparrow\mathbf{I}^{n-i}\colon \uparrow\mathbf{I}^{n+1} \to \uparrow\mathbf{I}^n.$

Now, the *directed cubical set* of a d-space $X$ and its left adjoint functor, the *directed geometric realisation* of a cubical set $K$, can be constructed as in classical case (cf. [Ma])

(3) $\quad \uparrow\mathcal{R}\colon \mathbf{Cub} \rightleftarrows \mathbf{dTop}\colon \uparrow\mathcal{C},$

$\quad\ \uparrow\mathcal{C}_n(X) = \mathbf{dTop}(\uparrow\mathbf{I}^n, X), \qquad\qquad \uparrow\mathcal{R}(K) = \int^{[n]} K_n\bullet\uparrow\mathbf{I}^n,$

the d-space $\uparrow\mathcal{R}(K)$ being the pasting in **dTop** of $K_n$ copies of $\uparrow\mathbf{I}^n$ ($n \geq 0$), along faces and degeneracies (precisely, the coend of the functor $K\bullet\uparrow\Delta\colon \mathbb{I}^{op}\times\mathbb{I} \to \mathbf{dTop}$).

The adjunction $U \dashv C^0$ (1.1) between spaces and d-spaces gives back the ordinary realisation $\mathcal{R} = U.\uparrow\mathcal{R}\colon \mathbf{Cub} \to \mathbf{Top}$, left adjoint to the ordinary cubical functor $\mathcal{C} = \uparrow\mathcal{C}.C^0\colon \mathbf{Top} \to \mathbf{Cub}$.

Various models of **dTop** are directed realisations of cubical sets: for the directed interval $\uparrow\mathbf{I}$, take $\uparrow\mathbf{i} = \{0 \to 1\}$; for $\uparrow\mathbf{R}$, take $\uparrow\mathbf{Z}$ with non degenerate 1-cubes $k \to k+1$; for $\uparrow\mathbf{O}^1$, take $\uparrow\mathbf{o}^1 = \{0 \rightrightarrows 1\}$; for $\uparrow\mathbf{S}^1$, take $\uparrow\mathbf{s}^1 = \{* \to *\}$.

**4.6. Local preorders and colimits.** We show now that such pastings cannot be realised within lp-spaces, and that colimits there (e.g. coequalisers) can fail.

The first step is to single out a geometric realisation $\uparrow\mathcal{R}(K)$ which is not of local preorder type. Consider the cubical set $K$ represented below (with two vertices $x, y$; two non-degenerate 1-cubes $a, b$; one non-degenerate 2-cube $w$)

(1) 

$\qquad\qquad\qquad K \qquad\qquad\qquad \uparrow C^-(\uparrow\mathbf{S}^1) \qquad\qquad\qquad T$

Its ordinary realisation is the compact disc $D = D(0, 1)$; its directed realisation $\uparrow\mathcal{R}(K)$ is the disc with the d-structure $\uparrow C^-(\uparrow\mathbf{S}^1)$ considered in 1.6.3; it can be viewed as the coequaliser, in **dTop**, of two edges $b_i\colon \uparrow\mathbf{I} \to T$ where the triangle $T$ has the structure $\uparrow\mathbf{I}^2/A$ obtained by collapsing the left edge $A = \partial_1^-(\mathbf{I})$ of the square $\uparrow\mathbf{I}^2$ (in **dTop** and lp**Top**).

Now, let us prove that $b_1$ and $b_2$ have no coequaliser in lp**Top**, and hence $K$ has no geometric realisation there. Suppose by absurd this coequaliser exists: it must be the disc $D = (D, \prec)$, with a suitable local preorder (because the forgetful functor lp**Top** $\to$ **Top** preserves all the existing colimits, having an obvious right adjoint, the chaotic preorder). $D$ inherits from $\prec$ a d-structure containing the universal one, $\uparrow C^-(\uparrow S^1)$ (because the projection $\uparrow I^2 \to D$ must factor through the former). Thus, as in 1.6, our relation $\prec$ must be chaotic on all concentric circles of some small disc $D(0, \varepsilon)$. Construct now a new lp-space $D' = (D, \prec')$ by a precedence relation which is chaotic on the disc $D(0, \varepsilon/2)$ and agrees with $\uparrow C^-(\uparrow S^1)$ on the complement. The projection $\Delta^2 \to D'$ is locally increasing, whence also the identity mapping $D \to D'$ must be so. But this is false around any point of $D(0, \varepsilon) \setminus D(0, \varepsilon/2)$.

Note the role played in this counterexample by the degenerate cube $e(x)$; and indeed, *pre*-cubical sets (without degeneracies) *can* be realised as locally ordered spaces, as proved in [FGR].

**4.7. Metrisability.** Directed spaces can be defined by 'asymmetric distances'. A generalised metric space $X$ in the sense of Lawvere [La], called here a *directed metric space* or *d-metric space*, is a set $X$ equipped with a *d-metric* $\delta: X \times X \to [0, \infty]$, satisfying the axioms

(1) $\quad \delta(x, x) = 0, \qquad\qquad\qquad\qquad \delta(x, y) + \delta(y, z) \geq \delta(x, z)$.

This structure is natural within the theory of enriched categories, as showed in [La]. (If the value $\infty$ is forbidden, $\delta$ is usually called a *quasi-pseudo-metric*, cf. [Ke]; including it has various structural advantages, e.g. the existence of all limits and colimits.)

d**Mtr** (or $\uparrow$**Mtr**) will denote the category of such d-metric spaces, with *d-contractions* f: $X \to Y$ ($\delta(x, x') \geq \delta(f(x), f(x'))$). Limits and colimits exist and are calculated as in **Set**; products have the $l_\infty$ d-metric and equalisers the restricted one, while sums have the obvious d-metric and coequalisers have the d-metric induced on the quotient:

(2) $\quad \prod_i X_i, \qquad \delta(x, y) = \sup_i \delta_i(x_i, y_i)$,

$\qquad \sum_i X_i, \qquad \delta((x, i), (y, i)) = \delta_i(x, y), \qquad \delta((x, i), (y, j)) = \infty \quad (i \neq j)$,

$\qquad X/R, \qquad \delta(\xi, \eta) = \inf_{\mathbf{x}}(\sum_i \delta(x_{2i-1}, x_{2i})) \qquad\qquad (x_1 \in \xi; \; x_{2i} R x_{2i+1}; \; x_{2n} \in \eta)$.

The *reflected* d-metric space $R(X) = X^{op}$ has the opposite d-metric, $\delta^{op}(x, y) = \delta(y, x)$. A *symmetric* d-metric ($\delta = \delta^{op}$), preferably written as d, is the same as an *écart* in [Bou].

A d-metric space $X = (X, \delta)$ has an associated bitopological space $(X, \tau^-, \tau^+)$. At the point $x_0 \in X$, the *past* topology $\tau^-$ (resp. the *future* topology $\tau^+$) has a canonical system of fundamental neighbourhoods consisting of *past* discs $D^-$ (resp. *future* discs $D^+$) centred at $x_0$

(3) $\quad D^-(x_0, \varepsilon) = \{x \in X \mid \delta(x, x_0) < \varepsilon\}, \qquad D^+(x_0, \varepsilon) = \{x \in X \mid \delta(x_0, x) < \varepsilon\}$.

This describes the forgetful functor to bitopological spaces, whence to d-spaces (via 1.4c)

(4) $\quad$ d**Mtr** $\to$ b**Top** $\to$ d**Top**,

$\qquad (X, \delta) \mapsto (X, \tau^-, \tau^+) \mapsto (X, \text{b}\mathbf{Top}(\uparrow I, X))$;

d-spaces which can be obtained in this way will be said to be *d-metrisable*.

The standard models of 1.2 are all d-metrisable, in a natural way. First, $\mathbf{R}^n$ will have the $l_\infty$-metric (written d), and $\mathbf{S}^1$ the geodetic one. For the directed real line $\uparrow\mathbf{R}$, take $\delta(x, x') = d(x, x')$ if $x \leq x'$, and $\infty$ otherwise; similarly for $\uparrow\mathbf{I}$, $\uparrow\mathbf{R}^n$, $\uparrow\mathbf{I}^n$, $\uparrow\mathbf{O}^1$. For the d-circle $\uparrow\mathbf{S}^1$, $\delta(x, x')$ is, as in 1.4, the length of the anticlockwise arc from $x$ to $x'$ (giving again the coequaliser of the maps 1.2.2-3, in d**Mtr**). Finally, $\uparrow\mathbf{S}^n$ can be realised as the coequaliser $\uparrow\mathbf{I}^n/\partial\mathbf{I}^n$ in the category d**Mtr**. Thus, for $\uparrow\mathbf{S}^2$, $\delta(x', x'')$ and $\delta(x'', x')$ are respectively the length of the solid and the dashed path

(5) 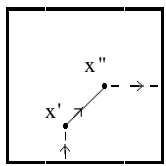